\documentclass[twocolumn]{autart}       

\usepackage{amsfonts}
\usepackage{amsmath}
\usepackage{epstopdf}
\usepackage{graphicx}

\newcommand{\norm}[1]{\lVert{#1}\rVert}

\newcommand{\R}{\mathbb{R}}

\newcommand{\N}{\mathbb{N}}

\newcommand{\ip}[2]{\left\langle #1, #2 \right\rangle}

\newcommand{\bmat}[1]{\begin{bmatrix}#1\end{bmatrix}}

\newcommand{\mcl}[1]{\mathcal{ #1}}
\newcommand{\mbf}[1]{\mathbf{ #1}}

\renewcommand{\norm}[1]{\left\lVert#1\right\rVert}

\newlength{\eqnspace}
\newlength{\parspace}
\begin{document}

		\begin{frontmatter}
	\title{A Partial Integral Equation Representation of Coupled Linear PDEs and Scalable Stability Analysis using LMIs}

	\author[Peet]{Matthew M. Peet}\ead{mpeet@asu.edu},               
	
	\address[Peet]{\textsc{School for the Engineering of Matter, Transport and Energy, Arizona State University, Tempe, AZ, 85298 USA.}}             
\thanks{This work was supported by grants NSF \# CNS-1739990, CMMI-1935453  and ONR \#N000014-17-1-2117.}


\begin{abstract}
We present a new Partial Integral Equation (PIE) representation of Partial Differential Equations (PDEs) in which it is possible to use convex optimization to perform stability analysis with little or no conservatism. The first result gives a standardized representation for coupled linear PDEs in a single spatial variable and shows that any such PDE, suitably well-posed, admits an equivalent PIE representation, defined by the given conversion formulae. This leads to a new prima facie representation of the dynamics without the implicit constraints on system state imposed by boundary conditions. The second result is to show that for systems in this PIE representation, convex optimization may be used to verify stability without discretization. The resulting algorithms are implemented in the Matlab toolbox PIETOOLS, tested on several illustrative examples, compared with previous results, and the code has been posted on Code Ocean. Scalability testing indicates the algorithm can analyze systems of up to 40 coupled PDEs on a desktop computer.  \vspace{-1mm}
\end{abstract}

	\end{frontmatter}

\section{Introduction}\vspace{-3mm}

Partial Differential Equations (PDEs) are used to model systems where the state depends continuously on both time and secondary independent variables. Common examples of such secondary dependence include space; as in flexible structures (Bernoulli-Euler beams) and fluid flow (Navier-Stokes); or maturation, as in cell populations and predator-prey dynamics.\vspace{\parspace}

The most common method for computational analysis of PDEs is to project the infinite-dimensional state onto a finite-dimensional vector space using, e.g. ~\cite{marion_1989,ravindran_2000,rowley_2005} and to use the extensive literature on control of Ordinary Differential Equations (ODEs) to test stability of, and design controllers for, the resulting finite-dimensional system. However, such discretization approaches are prone to instability (e.g. in the case of hyperbolic balance laws~\cite{karafyllis_book}), numerical ill-conditioning, and large-dimensional state-spaces. Furthermore, representation of PDEs using ODEs inevitably neglects higher-order modes, modes which can be inadvertently excited by feedback control via the well-known spillover effect~\cite{balas_1978}.\vspace{\parspace}

Work on \textit{computational} methods for analysis and control of PDEs which do not rely on discretization has been more limited. Perhaps the most well-known computational method for stabilization of PDEs without discretization is the backstepping approach to controller synthesis~\cite{krstic_book,smyshlyaev_2005,aamo_2013}. This approach is not optimization-based, however, and not typically used for stability analysis (An exception being~\cite{saba_2019}). Recently, there has been some work on the use of Linear Matrix Inequalities (LMIs) to find Lyapunov functions for linear and nonlinear PDEs - See~\cite{fridman_2009,fridman_2016,solomon_2015,gaye_2013}. However, because most of these works focus on the nonlinear case, the Lyapunov functions proposed therein are relatively simple and the resulting stability conditions conservative. An extension of the IQC framework to PDEs can be found in~\cite{matthieu_2020}.\vspace{\parspace}

Numerous analytic (non-computational) methods have been proposed over the years for analysis of PDEs, including the well-developed literature on Semigroup theory~\cite{lasiecka_book,curtain_book,bensoussan_book,karafyllis_book,bastin_book,luo_book} and the literature on Port-Hamiltonian systems~\cite{villegas_thesis} for selecting boundary inputs. However, these methods are typically ad-hoc - relying on the expertise of the user to propose and test energy metrics.\vspace{\parspace}

Recently, Sum-of-Squares (SOS) has been used for analysis and control of PDEs and examples can be found in~\cite{safi_2017,gahlawat_2017TAC,gahlawat_2016ACC,gahlawat_2015CDC} and~\cite{ahmadi_2016,valmorbida_2014,valmorbida_2016,gahlawat_2017}. While these SOS-based works are relatively accurate, they: 1) Are mostly limited to scalar PDEs; 2) Suffer from high computational complexity; 3) Are mostly ad-hoc, requiring significant effort to extend the results to new PDEs. For example, these methods have never been able to analyze stability of beam or wave equations. The source of the difficulty in using LMIs and SOS for stability analysis of PDEs is that the solution of a PDE is required to satisfy three sets of constraints: the differential equation; the boundary conditions; and continuity constraints. This is in contrast to ODEs, which are defined by bounded linear operators (matrices) and solutions to which need only satisfy a single differential equation.\vspace{\parspace}

\noindent\textbf{The Goal of the Paper} is to create, for PDEs, an equivalent of the LMI framework developed for ODEs. Historically, PDEs (as old as Newton) are defined by two conflicting sets of equations: the PDE itself, which moves the state; and the Boundary Conditions (BCs), which implicitly constrain the motion of the state. We want to unify these conflicting constraints into a new state-space representation of PDEs, defined by an algebra of bounded linear operators, and which directly incorporates: the PDE, the BCs, and the continuity constraints -- thereby eliminating issues of well-posedness and obviating the need to account for implicit constraints on the state.\vspace{\parspace}

This approach is fundamentally different than previous work using SOS or LMI-based methods. These previous efforts used SOS or positive matrices to propose candidate Lyapunov functions and then attempted to integrate the effect of boundary conditions into the derivative using, e.g. integration by parts. By contrast, our approach is to integrate the effect of boundary conditions into the \textit{dynamics} - thereby obviating the need to account for them in the stability analysis. As a result, our algorithms have no obvious source of conservatism and scale to systems of up to 40 coupled PDEs.\vspace{\parspace}


\noindent\textbf{\underline{Approach:}} In this paper, we propose the Partial Integral Equation (PIE) representation of PDEs. PIEs are infinite-dimensional state-space systems of the form\vspace{\eqnspace}
\begin{align}
\mcl T \dot{\mbf x}(t)&=\mcl A\mbf x(t),\notag\\
\mbf x(0)&=\mbf x_0 \in L_2^n[a,b],\\[-8mm]\label{eqn:PIE_0}
\end{align}
where the state, $\mbf x(t)$ is in $L_2^n[a,b]$, and the system parameters ($\mcl T, \mcl A$) are Partial Integral (3-PI) operators. 3-PI refers to the 3 matrix-valued parameters, $\{N_0,N_1,N_2\}$ which define every such operator $\mcl{P}_{\{N_i\}}:L_2[a,b]\rightarrow L_2[a,b]$ as\vspace{\eqnspace}
\begin{align*}
&\left(\mcl P_{\{N_0,N_1,N_2\}}\mbf x\right)(s):= N_0(s) x(s) ds\\
&\qquad \qquad   +\int_a^s N_1(s,\theta)x(\theta)d \theta+\int_s^bN_2(s,\theta)x(\theta)d \theta.\vspace{\eqnspace}
\end{align*}
As shown in Section~\ref{sec:algebra}, all 3-PI operators are $L_2$-bounded and together, they form an algebra (closed under addition, composition, scalar multiplication). Because they are algebraic, 3-PI operators inherit many of the properties of matrices and there is now a Matlab toolbox, PIETOOLS (using SOSTOOLS as a model), which allows for manipulation of 3-PI operators using matrix syntax and which can solve Linear PI Inequality (LPI) constrained optimization problems using the YALMIP syntax for LMIs.\vspace{\parspace}

\noindent\textbf{\underline{Converting PDE state to PIE state:}}
The first contribution of the paper (extending the results in~\cite{peet_2018CDC}) is to show that the solutions to a broad class of PDEs can be represented using PIEs. For this result, we propose an alternative state-space. To explain this change of state, we consider the following standardized representation of coupled linear PDEs in a single spatial variable, presented in Section~\ref{sec:PDE}.\vspace{\eqnspace}
\begin{align*}
\bmat{\dot{x}_0(t,s)\\\dot{x}_1(t,s)\\\dot{x}_2(t,s)} &= A_0(s)\underbrace{\bmat{x_0(t,s)\\x_1(t,s)\\x_2(t,s)}}_{\mbf x \in X}+A_1(s)\bmat{x_1(t,s)\\x_2(t,s)}_s\\[-3mm]
&\qquad \qquad \qquad \qquad +A_2(s)\bmat{x_2(t,s)}_{ss}\\[-8mm]\vspace{\eqnspace}
\end{align*}
with associated state-space\vspace{\eqnspace}\vspace{\eqnspace}
\begin{equation*}
X=\left\{\bmat{x_0\\x_{1}\\x_{2}}\in L_2^{n_0}\times H_1^{n_1} \times H_2^{n_2}\;:\; B {\scriptsize \bmat{ x_1(a) \\ x_1(b) \\ x_2(a) \\ x_2(b)\\ x_{2s}(a) \\ x_{2s}(b)}}=0\right\}.\vspace{\eqnspace}
\end{equation*}
Most 1D PDEs can be formulated using this standardized representation.\vspace{\parspace}

\noindent\textit{For example,} if we consider the damped wave equation\vspace{\eqnspace}
\begin{align*}
 \ddot u(t,s) &= u_{ss}(t,s)\hspace{-.5mm}-\hspace{-.5mm}2a \dot u(t,s)\hspace{-.5mm}-\hspace{-.5mm}a^2u(t,s), \, s \in [0,1] \\
  u(t,0) & =0,\qquad u_s(t,1)=-k \dot u(t,1)\vspace{\eqnspace}
\end{align*}
Then, setting $u_1=\dot u$ and $u_2=u$, we have\vspace{\eqnspace}
\[
\bmat{\dot u_1(t,s)\\ \dot u_2(t,s)}=\underbrace{\bmat{-2a&-a^2\\1&0}}_{A_0}\bmat{u_1(t,s)\\u_2(t,s)}+\underbrace{\bmat{1\\0}}_{A_2}u_{2ss}(t,s)\vspace{\eqnspace}
\]
with BCs $u_2(0)=0$, $u_1(0)=0$, and $u_{2s}(1)=-ku_1(1)$ so that\vspace{\eqnspace}
\begin{equation*}
X=\left\{\scriptsize\bmat{u_{1}\\u_{2}}\in H_1^{1} \times H_2^{1}: \underbrace{\scriptsize\bmat{0&0&1&0&0&0\\
      1&0&0&0&0&0\\
      0&k&0&0&0&1}}_{B}{\scriptsize\bmat{u_1(0)\\u_1(1)\\u_2(0)\\u_2(1)\\u_{2s}(0)\\u_{2s}(1)}=0}\right\}.\vspace{\eqnspace}
\end{equation*}\vspace{\parspace}

\noindent\textbf{\underline{The Fundamental State-Space}} A reasonable definition of state is the minimal amount of information needed to forward-propagate the solution. By this measure, and referring to the example, defining the state of a PDE as $u(t) \in H_1 \times H_2$ is not minimal, as this function contains redundant information regarding the boundary values. We propose, then, that for a PDE defined by $A_i$ and $X$, the correct definition of state is the so-called \textit{fundamental state}, where for $\mbf x\in X$, we define\vspace{\eqnspace}
\[
\mbf x_f=\bmat{x_0\\x_{1s}\\x_{2ss}}=\bmat{I&&\\&\partial_s&\\&&\partial_s^2}\mbf x\in L_2^{n_0+n_1+n_2}.\vspace{\eqnspace}
\]
In Section~\ref{sec:T}, we show that if the PDE is suitably well-posed, there exists a unitary 3-PI operator $\mcl T:L_2^{n_0+n_1+n_2}\rightarrow X$ such that $\mbf x=\mcl T \mbf x_f$.\vspace{\parspace}

\noindent\textbf{\underline{Equivalence of PIE and PDE:}} In Section~\ref{sec:PDE2PIE}, equipped with the unitary operator, $\mcl T$, we propose a 3-PI operator $\mcl A$ such that for any solution $\mbf x_f(t)\in L_2$ of the PIE in Eqn.~\eqref{eqn:PIE_0}, $\mbf x(t)=\mcl T\mbf x_f(t)$ satisfies the PDE defined by $A_i,X$ and that, conversely, for any solution of the PDE, $\mbf x$, the fundamental state $\mbf x_f(t) = \text{diag}(I, \partial_s,\partial_s^2)\mbf x(t)$ satisfies the PIE. We further show that exponential stability of the PIE and PDE are equivalent.\vspace{\parspace}

\noindent\textbf{\underline{An Linear PI Inequality (LPI) for Stability:}} Aside from have a minimal state-space, the advantage of the PIE representation of a PDE is computational. Recall that our goal is to create a framework akin to LMIs which can be used to study PDEs. Because PIEs do not have BCs, all information needed to define the solution is contained in the fundamental state, $\mbf x_f \in L_2$. Therefore, unlike PDEs, where the effect of the BCs needs to be ``brought in'' using integration by parts, Poincare inequalities, etc., PIEs are a prima facie representation of the dynamics. This allows us to pose the stability test for PIEs as a convex optimization problem of the following form:\vspace{\eqnspace}
\begin{align*}
&\substack{\text{Find}}_{\mcl P=\mcl P_{\{N_i\}}}:  \mcl P\ge \epsilon I,\\
&\mcl T^*\mcl P\mcl A+\mcl A^*\mcl P\mcl T\le -\delta \mcl T^*\mcl T.\\[-7mm]
\end{align*}
This LPI, then, is a straightforward generalization of the Lyapunov LMI: $P>0$, $A^TP+PA<0$. Furthermore, again motivated by the efficient Matlab parser YALMIP~\cite{yalmip}, we have recently developed the Matlab toolbox PIETOOLS allows us to manipulate 3-PI operators using matrix syntax and solve Linear PI Inequality (LPI) constrained optimization problems. Thus, if $\{\mcl T, \mcl A\}$ are as defined in Section~\ref{sec:PDE2PIE}, then our stability test reduces to 3 lines of Matlab code\vspace{\eqnspace}
\begin{flalign*}\scriptsize
&\texttt{[X,P] = sos\_posopvar(X,n,I,s,th,d);}&\\
&\texttt{D = -del*T'*T-A'*P*T-T'*P*A}&\\
&\texttt{X = sosopineq(X,D);}&\vspace{\eqnspace}
\end{flalign*}
where the functions $\texttt{sos\_posopvar}$ and $\texttt{sos\_opineq}$ are defined in Section~\ref{sec:PIETOOLS} and in~\cite{shivakumar_2020ACC}. We emphasize that this code applies to any PDE in standardized format and since there is no need to bring in the boundary conditions, there is no obvious source of conservatism. Specifically, in Section~\ref{sec:examples2}, we apply the algorithm to beam and wave equations for which there are no previous LMI-based stability conditions. Furthermore, the lack of conservatism is verified in Section~\ref{sec:examples1} by comparing against known stability margins taken from the literature.\vspace{\parspace}

Finally, computational complexity depends on the degree of the polynomial parameters in the 3-PI variable $\mcl P$ (corresponding to the complexity of the candidate Lyapunov funtion). However, most problems only require very simple Lyapunov functions. In this case, the scalability of the method is comparable to the complexity of discretization-based analysis. Specifically, if we choose the polynomial parameters to have degree 2, then the algorithm can analyze stability of more than 40 coupled PDEs on a desktop computer.\vspace{\parspace}

In the following two illustrations, we attempt to further motivate our BC free PIE representation by illustrating how BCs affect the dynamics and complicate stability analysis in the original PDE state.\vspace{-2mm}

\subsection{BCs Complicate Stability Analysis}\vspace{-3mm}
The goal of the paper is to find a representation of PDEs which can be interpreted as we would interpret an ODE - without implicit constraints on the solution imposed by BCs or continuity. To motivate this goal, and to illustrate that BCs cannot be ignored in the PDE framework, let us illustrate what happens when we treat a PDE like an ODE - i.e. \textit{without} considering continuity constraints or boundary conditions. Obviously, in a more classical analysis, we may account for such boundary conditions using, e.g. integration-by-parts. However, remember that our goal in this paper is to eliminate the need for such ad hoc manipulations - as is possible in the PIE representation. To this end, suppose we are given a PDE, parameterized by $A_0,A_1,A_2$ as follows.\vspace{\eqnspace}
\[
\dot{\mbf x}(t,s)=A_0(s)\mbf x(t,s)+A_1(s)\mbf x_s(t,s) +A_2(s)\mbf x_{ss}(t,s)\vspace{\eqnspace}
\]
An obvious class of candidate Lyapunov functions for this system is parameterized by $M$ as\vspace{\eqnspace}
 \[
V(\mbf x)=\ip{\mbf x}{\mcl P_{\{M,0,0\}}\mbf x}_{L_2}=\int_{a}^b \mbf x(s)^TM(s)\mbf x(s)ds.\vspace{\eqnspace}
\]
As would be the case for an ODE, $V(\mbf x)\ge \epsilon \norm{\mbf x}^2$ if $M(s)\ge \epsilon I$ for all $s$ and some $\epsilon >0$ - a constraint which is easy to enforce using, e.g. SOS. However, if we now take the derivative of this candidate Lyapunov function we obtain\vspace{\eqnspace}\vspace{-1mm}
\begin{align*}
&\dot V(\mbf x)=\int\limits_{a}^b \bmat{\mbf x(s)\\\mbf x_s(s)\\\mbf x_{ss}(s)}^TD(s)\bmat{\mbf x(s)\\\mbf x_s(s)\\\mbf x_{ss}(s)}ds\\
&D(s):={\scriptsize\bmat{A_0(s)^TM(s)+M(s)A_0(s)&M(s)A_1(s)&M(s)A_2(s)\\A_1(s)^TM(s)&0&0\\A_2(s)^TM(s)&0&0}}.\vspace{\eqnspace}
\end{align*}
Now, if we were to treat the system like an ODE (without considering the BCs and continuity constraints), we would constrain $D(s)\le 0$ and this would imply stability. Unfortunately, however, $D(s)\not\le0$ for ANY choice of $M,A_1,A_2\neq 0$. The problem, of course, is that the differentiation operator branches $\mbf x$ into $\mbf x_s$ and $\mbf x_{ss}$, neither of which are independent of $\mbf x$. Moreover, the information which determines the relationship between $\mbf x$, $\mbf x_s$ and $\mbf x_{ss}$ is not embedded in the differential equation. Rather, this information is implicit in the BCs and continuity constraints. At this point, of course, the user would attempt to  ``bring in'' the boundary and continuity properties to obtain a new stability condition using, e.g. integration by parts or Stokes Theorem~\cite{gahlawat_2017}. However, such secondary steps require significant expertise, are ad-hoc, and may result in conservatism. In the analysis of ODEs, by contrast, no such secondary analysis is necessary. The goal of this paper, then, is to find a prima facie representation of the PDE with no implicit constraints on continuity and BCs, and wherein we may develop a computational framework which mirrors that of LMIs for ODEs.\vspace{-2mm}

\subsection{BCs Significantly alter the Dynamics}\vspace{-3mm}
Having observed that BCs complicate the analysis of PDEs, let us now examine how the BCs can alter the nature of a distributed parameter system. Specifically, in this subsection, we consider an example of how the incorporation of continuity and BCs can take a system (wherein the dynamics do not explicitly depend on the partial derivatives), and transform it into a system with explicit dependence on the partial derivatives. Consider the following non-partial-differential, yet distributed-parameter system.\vspace{\eqnspace}
\[
\dot{\mbf u}(t,s)= \mbf u(t,s),\; \mbf u(t,0)=w_1(t),\; \mbf u_s(t,0)=w_2(t).\vspace{\eqnspace}
\]
To allow for the specified BCs, we restrict continuity of $\mbf u$ as $\mbf u(t) \in H_2^1$. Note that this system is not, prima facie, a PDE and that the dynamics are identical at every point in the domain. However, if we now combine the fundamental theorem of calculus with integration by parts, we obtain a very different, yet equivalent, set of dynamics - dynamics with explicit dependence on the partial derivatives.\vspace{\eqnspace}\vspace{-1mm}
\[
\dot{\mbf u}(t,s)=sw_1(t)+w_2(t)+\int^s_0 (s-\eta) \mbf u_{\eta \eta}(t,\eta)d \eta\vspace{\eqnspace}
\]
This formulation of the same system directly incorporates continuity and BCs into the dynamics - which are now expressed using the partial derivative $\mbf u_{ss}$. In addition, while the original formulation was spatially decoupled, with $w_1,w_2$ only acting at the boundary, the new formulation shows that the exogenous inputs $w_1,w_2$ are felt instantaneously at every point in the domain. This observation demonstrates that the continuity properties of the solution (combined with BCs) fundamentally alter the dynamics.\vspace{\parspace}

The first goal of this paper, then, is: 1) to show that a broad class of PDEs can be reformulated in a way which specifies precisely how the BCs affect the dynamics and 2) to provide universal formulae for constructing such a representation.
\vspace{-3mm}

\section{Notation}\vspace{-3mm}
We define $L_2[a,b]^n$ to be space of $\R^n$-valued Lesbegue integrable functions defined on $[a,b]$ and equipped with the standard inner product.
We use $W^{k,p}[a,b]^n$ to denote the Sobolev subspace of $L_p[a,b]^n$ defined as $\{u \in L_p[a,b]^n\, :\, \frac{\partial^q}{\partial x^q} u \in L_{p} \text{ for all } q \le k \}$. $H_k[a,b]:=W^{k,2}[a,b]$ and $H_k^n[a,b]=\prod_{i=1}^n H_k[a,b] $. For efficiency, we typically omit the domain, so that, e.g. $H_k^n:=H_k^n[a,b]$ unless otherwise stated. $I_n\in \R^{n \times n}$ and $0_{n_1 \times n_2}\in \R^{n_1\times n_2}$ are used to denote the identity and zero matrices and the subscripts are omitted when the dimension of the matrices is clear from context. $\mbf I$ denotes the indicator function $\mbf I:\R \rightarrow \{0,1\}$, defined as\vspace{\eqnspace}
\[
\mbf I(s)=\begin{cases}
       1, & \mbox{if } s>0 \\
       0, & \mbox{otherwise}.
     \end{cases}\vspace{\eqnspace}\vspace{-3mm}
\]

\section{A Standardized PDE Representation}\label{sec:prelims}\label{sec:PDE} \vspace{-3mm}
The two primary contributions of this paper are: a formula for conversion of PDEs to PIEs; and an LPI framework for Lyapunov stability analysis of PIEs (Section~\ref{sec:LOI}). The significance of the latter result clearly depends on the set of PDEs which can be converted to PIEs. In this section, we propose a standardized framework for representation of PDEs. In Section~\ref{sec:PDE2PIE}, we will show that for any such PDE, there exists a PIE for which a solution to the PIE yields a solution to the PDE and vice-versa. The class of PDEs considered here is not exhaustive, however. That is, there exist PDEs not listed here which may be converted to PIEs. Furthermore, there exist PIEs which do not have a coupled PDE representation.\vspace{\parspace}

We consider coupled PDEs with infinitesimal generator of the following form\vspace{\eqnspace}
\begin{align}
\bmat{\dot{x}_0(t,s)\\\dot{x}_1(t,s)\\\dot{x}_2(t,s)} &= A_0(s)\bmat{x_0(t,s)\\x_1(t,s)\\x_2(t,s)}+A_1(s)\bmat{\partial_s x_1(t,s)\\ \partial_s x_2(t,s)}\notag \\
&\qquad \qquad \qquad +A_2(s)\bmat{\partial_s^2 x_2(t,s)}. \label{eqn:PDE_1}\\[-7mm]\notag
\end{align}
and with domain\vspace{\eqnspace}
\begin{equation}
X:=\left\{\bmat{x_0\\x_{1}\\x_{2}}\in L_2^{n_0}\times H_1^{n_1} \times H_2^{n_2}\;:\; B {\scriptsize \bmat{ x_1(a) \\ x_1(b) \\ x_2(a) \\ x_2(b)\\ x_{2s}(a) \\ x_{2s}(b)}}=0\right\}\label{eqn:X_def}\vspace{\eqnspace}
\end{equation}
where\vspace{\eqnspace}
\begin{equation}
B\bmat{I_{n_1}&0 &0\\I_{n_1}&0&0 \\0&I_{n_2} &0\\0&I_{n_2}&(b-a)I_{n_2}\\0&0&I_{n_2}\\0&0 &I_{n_2}}\qquad \text{is invertible.}\label{eqn:B_assumption}\vspace{\eqnspace}
\end{equation}
Specifically, given $\mbf x_0 \in X$, we say that $\mbf x$ satisfies the PDE defined by $\{A_i,X\}$ if $\mbf x$ is Frech\'et differentiable, $\mbf x(0)=\mbf x_0$, $\mbf x(t)\in X$ and Equation~\eqref{eqn:PDE_1} is satisfied for all $t\ge 0$.\vspace{-2mm}

\subsection{A Guide to Partition of States}\vspace{-3mm}
The partition of states in Equation~\eqref{eqn:PDE_1} is not overly restrictive and there is no special structure to this generator. The partition is purely organizational and defined by the domain restriction $\mbf x \in L_2^{n_0}\times H_1^{n_1} \times H_2^{n_2}$. States with no restrictions on continuity are assigned to be an element of the vector $x_0(t)\in L_2[a,b]^{n_0}$. If a state has no continuity properties, then these states cannot be differentiated and it is not possible to assign boundary conditions, as the limits $x_0(a), x_0(b)$ do not exist. States which are continuous, but not continuously differentiable are assigned to $x_1(t)\in H_{1}^{n_1}$. The continuity property of these states allow for boundary conditions, as $x_1(a), x_1(b)$ exist. However, since $\partial_s x_1\in L_2[a,b]^{n_1}$ is not continuous, we cannot assign boundary conditions which involve $x_{1s}(a)$ or $x_{1s}(b)$, as these limits do not exist. Finally, states which are required to be continuously differentiable are assigned to $x_2(t)\in H_{2}^{n_2}$ and admit boundary conditions involving $x_{2s}(a)$ or $x_{2s}(b)$ and second-order spatial derivatives, $\partial^2_s x_{2}$. This standardized representation specifically excludes states in $H^n_k$ where $k>2$. Although the results of the paper can be extended to such systems, such an extension is not considered here.\vspace{-2mm}
\subsection{A Guide to Boundary Conditions}\vspace{-3mm}
In this subsection, we propose restrictions on the matrix $B$ which are equivalent to Equation~\eqref{eqn:B_assumption}. Specifically, the row rank of $B$ must be $n_1+2n_2$ and $B$ contains no boundary conditions of a given prohibited form. Note that the rank condition on $B$ is not overly restrictive as, to the best of our knowledge, it is a necessary condition for existence of a unique solution for any PDE in standardized form.\vspace{-2mm}

\subsubsection{Prohibited Boundary Conditions}\vspace{-3mm}
A necessary and sufficient condition for $B$ to satisfy Equation~\eqref{eqn:B_assumption} is for $B\in \R^{(n_1+2n_2)\times (2n_1+4n_2)}$ to have row rank $n_1+2n_2$ and to define no boundary conditions consisting of a linear combination of $x_1(a)-x_1(b)=0$, $x_2(a)+(b-a)x_{2s}(a)-x_2(b)=0$, or $x_{2s}(a)-x_{2s}(b)=0$.

\begin{lem} Suppose $B\in \R^{(n_1+2n_2)\times (2n_1+4n_2)}$. Define\vspace{\eqnspace}
\[
T^\perp:=\bmat{I_{n_1}&-I_{n_1}&0&0&0&0\\0&0&I_{n_2}&-I_{n_2}&I_{n_2}(b-a)&0\\0&0&0&0&I_{n_2}&-I_{n_2}}.\vspace{\eqnspace}
\]
Equation~\eqref{eqn:B_assumption} is satisfied if and only if $B\in \R^{(n_1+2n_2)\times (2n_1+4n_2)}$,  has row rank $n_1+2n_2$ and the row space of $B$ and the row space of $T^\perp$ has trivial intersection.\vspace{-3mm}
\end{lem}
\begin{pf}
Let the $T$ matrix from Eqn.~\eqref{eqn:B_assumption} be as defined in Eqn.~\eqref{eqn:Gdefs}.
Now suppose $BT$ is invertible. Since $T \in \R^{(2n_1+4n_2) \times (n_1+2n_2)}$, we require $B\in \R^{(n_1+2n_2)\times (2n_1+4n_2)}$ in order for $BT$ to exist and be square. Since $BT \in \R^{(n_1+2n_2) \times (n_1+2n_2)}$, $B$ must also have row rank $n_1+2n_2$. Now, since $T$ has column rank $n_1+2n_2$, its row rank is also $n_1+2n_2$. Now $T^\perp$ has row rank $n_1+2n_2$ and $T^\perp T=0$. Thus the row space of $T^\perp$ lies in $Im(T)^\perp$ and since $Im(T)^\perp$ is of dimension $n_1+2n_2$, the row space of $T^\perp$ is $Im(T)^\perp$. Therefore $x^TBT=0$ for some $x\neq0$, if and only if the intersection of the row space of $B$ and that of $T^\perp$ is non-trivial. This establishes necessity. For sufficiency, we assume $B\in \R^{(n_1+2n_2)\times (2n_1+4n_2)}$ and has row rank $n_1+2n_2$. Furthermore, as shown, the row space of $B$ has trivial intersection with $Im(T)^\perp$. Again, we have that $x^TBT=0$ implies $x=0$, from which we conclude invertibility. $\square$\vspace{-3mm}
\end{pf}
\begin{note}
The restriction on prohibited boundary condition is subtle. For example, $x_2(a)=x_2(b)$ is permitted, except if combined with $x_{2s}(a)=0$ (combining with $x_{2s}(b)=0$ is still OK). Meanwhile, $x_1(a)=x_1(b)$ and $x_{2s}(a)=x_{2s}(b)$ are never OK. Additionally, $x_{2s}(a)=0$ is OK, unless combined with $x_2(a)=x_2(b)$ or $x_{2s}(b)=0$. Of course, the most reliable way to check if certain boundary conditions are permitted is to simply construct $B$ and check the rank of $BT$. The PIETOOLS implementation described in Section~\ref{sec:PIETOOLS} will do this automatically and generate an error if $BT$ is not invertible.\vspace{-3mm}
\end{note}

\subsubsection{A Note on Necessity of Equation~\eqref{eqn:B_assumption}}\vspace{-3mm}
Boundary conditions of the form $x_1(a)=x_1(b)$ are periodic and imply\vspace{\eqnspace}\vspace{-2mm}
\[
\int_a^b x_{1s}(s)ds=0.\vspace{\eqnspace}
\]
Likewise $x_{2s}(a)=x_{2s}(b)$ implies\vspace{\eqnspace}
\[
\int_a^b x_{2ss}(s)ds=0,\vspace{\eqnspace}
\]
and $x_2(a)+(b-a)x_{2s}(a)=x_2(b)$ implies\vspace{\eqnspace}
\[
\int_a^b \int_{a}^s x_{2\eta \eta}(\eta)d \eta ds=0.\vspace{\eqnspace}
\]
In this way, the prohibited BCs represent integral constraints on the respective PIE (fundamental) states, $x_{1s} \in L_2$ and $x_{2ss}\in L_2$, meaning these PIE states are not minimal (dynamics expressed using these states have implicit constraints). One option in these cases may be to redefine the PIE states modulo an integral constraint, however this extension is left for future work.\vspace{-2mm}

\subsection{Euler-Bernoulli Beam Example}\label{subsec:EB}\vspace{-4mm}
In order to better understand how to write a PDE in the standardized PDE form of Eqn.~\eqref{eqn:PDE_1}, let us consider the cantilevered Euler-Bernoulli beam:\vspace{\eqnspace}
\begin{align*}
&\ddot{u}(t,s)=-cu_{ssss}(t,s),\qquad \text{where}\quad\\
& u(0)=u_{s}(0)=u_{ss}(L)=u_{sss}(L)=0.\\[-7mm]
\end{align*}
We wish to construct a standardized PDE representation of this classic diffusive model. Following the approach in, e.g.~\cite{villegas_thesis} (from which we also get the Timoshenko beam model in Section~\ref{sec:examples2}), we first introduce the augmented state $u_1=\dot{u}$ - a choice which leads to ``natural'' BCs for which the system is well-posed~\cite{luo_book}. This choice also eliminates the second order time-derivative, $\ddot{u}$. Since $u \in H_4$, we eliminate the fourth-order spatial derivative by creating the augmented state $u_2=u_{ss}$. Taking the time-derivative of these states, we obtain\vspace{\eqnspace}
\begin{align*}
  \dot u_1 & =\ddot{u}=-cu_{ssss}=-cu_{2ss} \\
  \dot u_{2} & =\partial_t \partial_s^2 u =\partial_s^2 \dot{u}=u_{1ss}.\\[-7mm]
\end{align*}
These equations are now in the standardized form\vspace{\eqnspace}
\[
\dot{\mbf x}(t)=\underbrace{\bmat{0&-c\\1&0}}_{A_2}\mbf x_{ss}(t)\vspace{\eqnspace}
\]
where $A_0=A_1=0$, $n_2=2$, and $n_0=n_1=0$. We now examine the boundary conditions using these states:\vspace{\eqnspace}
\[
u_{ss}(L)=u_2(L)=0\quad \text{and}\quad u_{sss}(L)=u_{2s}(L)=0.\vspace{\eqnspace}
\]
These boundary conditions are insufficient, as the resulting rank is 2. However, we may differentiate boundary conditions in time to obtain new boundary conditions\vspace{\eqnspace}
\[
\dot{u}(0)=u_1(0)=0\quad \text{and}\quad \dot{u}_{s}(0)=u_{1s}(0)=0.\vspace{\eqnspace}
\]
We now have 4 boundary conditions, which we use to construct the $B$ matrix as\vspace{\eqnspace}\vspace{-3mm}
\[
\underbrace{\bmat{1&0&0&0&0&0&0&0\\
0&0&0&1&0&0&0&0\\
0&0&0&0&1&0&0&0\\
0&0&0&0&0&0&0&1}}_{B}{\scriptsize\bmat{u_1(0)\\u_2(0)\\u_1(L)\\u_2(L)\\u_{1s}(0)\\u_{2s}(0)\\u_{1s}(L)\\u_{2s}(L)}}=0.\vspace{\eqnspace}
\]
The $B$ matrix is now of rank $4=n_{1}+2n_2$ and satisfies Equation~\eqref{eqn:B_assumption}.\vspace{\parspace}

Now, if $u$ satisfies the E-B beam equation for some initial condition, we have that $u_1=\dot u,u_2=u_{ss}$ satisfy the standardized PDE model. However, conversely, if $u_1, u_2$ satisfy the standardized PDE for some initial condition, then in order to construct a solution to the original PDE, we must integrate $u_1$ in time. However, this requires knowledge of $u(0)$. Thus we find that some information on the system solution has been lost in the standardized representation. Note, however, that we could retain this information by including a third state, $u_3=u$, so that $\dot u_3=u_1$ and then the solutions would be equivalent.\vspace{-2mm}

\subsection{Exponential Stability of Coupled PDE Systems}\vspace{-3mm}
In this subsection, we define two notions of exponential stability with respect to the standardized PDE representation - stability in the $X$ norm and stability in the $L_2$ norm. In Section~\ref{sec:PIE}, we will define the notion of exponential stability for PIEs. In Section~\ref{sec:LOI}, we will show that exponential stability of a PIE representation of a standardized PDE is equivalent to exponential stability of the original standardized PDE in the $X$ norm.\vspace{-2mm}

\begin{defn}
We say the PDE defined by $\{A_i,X\}$ is exponentially stable in $X$ if there exist constants $K,\gamma>0$ such that for any $\mbf x_0 \in X$, any solution $\mbf x$ of the PDE defined by $\{A_i,X\}$ satisfies\vspace{\eqnspace}
\[
\norm{\mbf x(t)}_{L_2^{n_0}\times H_1^{n_1} \times H_2^{n_2}}\le K\norm{\mbf x_0}_{L_2^{n_0}\times H_1^{n_1} \times H_2^{n_2}}e^{-\gamma t}.\vspace{\eqnspace}
\]
\end{defn}

\begin{defn}
We say the PDE defined by $\{A_i,X\}$ is exponentially stable in $L_2$ if there exist constants $K,\gamma>0$ such that for any $\mbf x_0 \in X$, any solution $\mbf x$ of the PDE defined by $\{A_i,X\}$ satisfies\vspace{\eqnspace}
\[
\norm{\mbf x(t)}_{L_2^{n_0+n_1+n_2}}\le K\norm{\mbf x_0}_{L_2^{n_0+n_1+n_2}}e^{-\gamma t}.\vspace{\eqnspace}
\]
\end{defn}

\begin{note}
Exponential stability in $X$ implies exponential stability in $L_2$, since $\norm{\mbf x}_{L_2}\le \norm{\mbf x}_{H_k}$ for any $\mbf x \in H_k$ and $k\ge 0$. Furthermore, our definitions of exponential stability implies that all states in $\mbf x$ must be exponentially decreasing in the given norm. Because not all standardized PDE representations of a given scalar high-order PDE necessarily use the same set of first-order states (See, e.g. the E-B beam example), this raises the possibility one standardized PIE representation may be exponentially stable, while another may not.\vspace{-2mm}
\end{note}

Note that in the case where the $X$ or $L_2$ stability definition holds with $\gamma=0$, we say that the system is Lyapunov stable or neutrally stable.\vspace{-2mm}

\section{3-PI Operators Form an Algebra}\label{sec:algebra}\vspace{-3mm}
In Section~\ref{sec:PDE2PIE}, we will construct a PIE representation of any PDE in the standardized form described in Section~\ref{sec:PDE}. PIEs, as will be defined in Section~\ref{sec:PIE}, have the advantage that they are parameterized by the class of 3-PI operators, which are bounded on $L_2$ and form an algebra. Furthermore, candidate Lyapunov functions can be parameterized using 3-PI operators. The algebraic nature of 3-PI operators significantly simplifies the problem of analysis and control of PIEs.\vspace{\parspace}

Formally, we say that an operator $\mcl P$ is 3-PI if
there exist 3 bounded matrix-valued functions $N_0:[a,b]\rightarrow \R^{n \times n}$, $N_1:[a,b]^2 \rightarrow \R^{n \times n}$, and $N_2:[a,b]^2 \rightarrow \R^{n \times n}$ such that\vspace{\eqnspace}
\begin{align*}
&(\mcl P \mbf x)(s):=\left(\mcl P_{\{N_0,N_1,N_2\}}\mbf x\right)(s):= N_0(s) x(s)  \\
& \qquad +\int_a^s N_1(s,\theta)x(\theta)d \theta+\int_s^bN_2(s,\theta)x(\theta)d \theta,\\[-7mm]
\end{align*}
where $N_0$ defines a multiplier operator and $N_1,N_2$ define the kernel of an integral operator.\vspace{\parspace}

For given $N_0, N_1,N_2$, we use $\mcl P_{\{N_0,N_1,N_2\}}:L_2^n \rightarrow L_2^n$ to denote the corresponding PI operator. When clear from context, we use the shorthand notation $\mcl P_{\{N_i\}}$ to indicate $\mcl P_{\{N_0,N_1,N_2\}}$.\vspace{\parspace}

One may interpret 3-PI operators to be an extension of matrices, wherein $N_0$ defines the diagonal of the matrix, $N_1$ contains the sub-diagonal terms, and $N_2$ contains the terms above the diagonal.\vspace{\parspace}

In the following subsections, we show that this class of bounded linear operators is closed under composition and adjoint (closure under addition and scalar multiplication follows immediately from addition and scalar multiplication of parameters). Furthermore, these results show that if we define the set of 3-PI operators with polynomial parameters $N_0$, $N_1$, and $N_2$, then this set forms a sub-algebra. \vspace{-4mm}

\subsection{Composition of 3-PI operators}\vspace{-3mm}
In this subsection, we derive an analytic formula for the composition 3-PI operators. Specifically, we have the following.\vspace{-3mm}

\begin{lem}\label{thm:composition}
 For any bounded functions $B_0,N_0:[a,b]\rightarrow \R^{n \times n}$, $B_1,B_2,N_1,N_2:[a,b]^2 \rightarrow \R^{n \times n}$, we have\vspace{\eqnspace}
\[
\mcl P_{\{R_i\}}=  \mcl P_{\{B_i\}}  \mcl P_{\{N_i\}}\vspace{\eqnspace}\vspace{-2mm}
\]
where\vspace{\eqnspace}
\begin{align}
&R_0(s)=B_0(s)N_0(s), \label{eqn:composition}\\
&R_1(s,\theta)=B_0(s)N_1(s,\theta)+B_1(s,\theta)N_0(\theta) \notag\\
&\; +\int_{a}^\theta B_1(s,\xi) N_2(\xi,\theta)d \xi
+\int_{\theta}^s B_1(s,\xi) N_1(\xi,\theta)d \xi\notag \\
&\qquad \qquad \qquad \qquad \qquad  + \int_{s}^b B_2(s,\xi) N_1(\xi,\theta)d \xi \notag, \\
&R_2(s,\theta)=B_0(s)N_2(s,\theta)+B_2(s,\theta)N_0(\theta) \notag \\
&\; +\int_{a}^s B_1(s,\xi) N_2(\xi,\theta)d \xi
+\int_{s}^\theta B_2(s,\xi) N_2(\xi,\theta)d \xi\notag \\
&\qquad \qquad \qquad \qquad \qquad \;+\int_{\theta}^b B_2(s,\xi) N_1(\xi,\theta)d \xi. \notag\\[-9mm]\notag
\end{align}
\end{lem}
\vspace{-3mm}
\begin{pf}
To prove the theorem, we exploit the linear structure of the operator to decompose\vspace{\eqnspace}
\[
  \mcl P_{\{B_0,B_1,B_2\}} =\mcl P_{\{B_0,0,0\}}+\mcl P_{\{0,B_1,0\}}+\mcl P_{\{0,0,B_2\}}\vspace{\eqnspace}
\]
and\vspace{\eqnspace}
\[
  \mcl P_{\{N_0,N_1,N_2\}} =\mcl P_{\{N_0,0,0\}}+\mcl P_{\{0,N_1,0\}}+\mcl P_{\{0,0,N_2\}}.\vspace{\eqnspace}
\]
Then\vspace{\eqnspace}
  \begin{align*}
  &\mcl P_{\{B_0,B_1,B_2\}}  \mcl P_{\{N_0,N_1,N_2\}}\\
  &=\mcl P_{\{B_0,0,0\}}\mcl P_{\{N_0,N_1,N_2\}}+\mcl P_{\{0,B_1,B_2\}}\mcl P_{\{N_0,0,0\}}\\
  &\quad +\mcl P_{\{0,B_1,0\}}  \mcl P_{\{0,N_1,0\}}+\mcl P_{\{0,B_1,0\}}  \mcl P_{\{0,0,N_2\}}\\
  &\quad +\mcl P_{\{0,0,B_2\}}  \mcl P_{\{0,N_1,0\}}+\mcl P_{\{0,0,B_2\}}  \mcl P_{\{0,0,N_2\}}.
  \end{align*}
We now consider each term separately, starting with the first two, which are trivial. First,\vspace{\eqnspace}
\begin{align*}
&\left(\mcl P_{\{B_0,0,0\}}  \mcl P_{\{N_0,N_1,N_2\}}x\right)(s)\\
&=B_0(s)N_0(s)x(s)+\int_{a}^s B_0(s)N_1(s,\theta)x(\theta)d \theta\\
&\qquad \qquad +\int_{s}^b B_0(s)N_2(s,\theta)x(\theta)d \theta\\
&=\mcl P_{\{R_{0},R_{1a},R_{2a}\}},\vspace{\eqnspace}
\end{align*}
where\vspace{\eqnspace}
\begin{align*}
&R_{0}(s)=B_0(s)N_0(s), \quad R_{1a}(s,\theta)=B_0(s)N_1(s,\theta),\\
&R_{2a}(s,\theta)=B_0(s)N_2(s,\theta).\vspace{\eqnspace}
\end{align*}
 Similarly,{
  \begin{align*}
&  \left(\mcl P_{\{0,B_1,B_2\}}  \mcl P_{\{N_0,0,0\}}x\right)(s)\\
&=\int_{a}^s B_1(s,\theta)N_0(\theta)x(\theta)d \theta+\int_{s}^b B_2(s,\theta)N_0(\theta)x(\theta)d \theta\\
&=\mcl P_{\{0,R_{1b},R_{2b}\}},
  \end{align*}}
where\vspace{\eqnspace}
\[
R_{1b}(s,\theta)=B_1(s,\theta)N_0(\theta),\; R_{2b}(s,\theta)=B_2(s,\theta)N_0(\theta).\vspace{\eqnspace}
\]
We now proceed to the more difficult terms. For these, recall the indicator function\vspace{\eqnspace}
\[
\mbf I(s)=\begin{cases}
       1, & \mbox{if } s>0 \\
       0, & \mbox{otherwise}.
     \end{cases}\vspace{\eqnspace}
\]
For the first term, we note that\vspace{\eqnspace}
\[
\mbf I(s-\eta) \mbf I(\eta-\xi) =
\begin{cases}
       \mbf I(s-\xi), & \mbox{if } \eta\in [\xi,s] \\
       0, & \mbox{otherwise}.\vspace{\eqnspace}
\end{cases}
\]
This allows us to change the variables of integration as follows.\vspace{\eqnspace}{
\begin{align*}
 & \left(\mcl P_{\{0,B_1,0\}}  \mcl P_{\{0,N_1,0\}}x\right)(s)\\
&=\int_{a}^s B_1(s,\eta)\int_{a}^\eta N_1(\eta,\xi)x(\xi) d \xi d \eta\\
&=\int_{a}^b \mbf I(s-\eta) B_1 (s,\eta) \int_{a}^b \mbf I(\eta-\xi) N_1(\eta,\xi)x(\xi) d \xi d \eta\\
&=\int_{a}^b \left(\int_{a}^b \mbf I(s-\eta) \mbf I(\eta-\xi) B_1(s,\eta)  N_1(\eta,\xi)d \eta \right) x(\xi) d \xi \\
&=\int_{a}^b \left(\int_{\xi}^s \mbf I(s-\xi) B_1(s,\eta) N_1(\eta,\xi)d \eta \right) x(\xi) d \xi \\
&=\int_{a}^s \left(\int_{\xi}^s  B_1(s,\eta) N_1(\eta,\xi)d \eta \right) x(\xi) d \xi \\
&=  \left(\mcl P_{\{0,R_{1c},0\}} x\right)(s),
  \end{align*}}
where\vspace{\eqnspace}
\[
R_{1c}(s,\theta)=\int_{\theta}^s B_1(s,\xi) N_1(\xi,\theta)d \xi.\vspace{\eqnspace}
\]
Next, we use another identity\vspace{\eqnspace}
\[
\mbf I(s-\eta) \mbf I(\xi-\eta) =\mbf I(s-\xi)\mbf I(\xi-\eta) + \mbf I(\xi-s)\mbf I(s-\eta)\vspace{\eqnspace}
\]
to establish the following.\vspace{\eqnspace}{\scriptsize
\begin{align*}
 & \left(\mcl P_{\{0,B_1,0\}}  \mcl P_{\{0,0,N_2\}}x\right)(s)\\
&=\int_{a}^s B_1(s,\eta)\int_{\eta}^b N_2(\eta,\xi)x(\xi) d \xi d \eta\\
&=\int_{a}^b \mbf I(s-\eta) B_1 (s,\eta) \int_{a}^b \mbf I(\xi-\eta) N_2(\eta,\xi)x(\xi) d \xi d \eta\\
&=\int_{a}^b \left(\int_{a}^b \mbf I(s-\eta) \mbf I(\xi-\eta) B_1(s,\eta)  N_2(\eta,\xi)d \eta \right) x(\xi) d \xi \\
&=\int_{a}^b \mbf I(s-\xi) \left(\int_{a}^\xi  B_1(s,\eta) N_2(\eta,\xi)d \eta \right) x(\xi) d \xi \\
&\qquad+ \int_{a}^b \mbf I(\xi-s) \left(\int_{a}^s  B_1(s,\eta) N_2(\eta,\xi)d \eta \right) x(\xi) d \xi \\
&=\int_{a}^s \left(\int_{a}^\xi  B_1(s,\eta) N_2(\eta,\xi)d \eta \right) x(\xi) d \xi \\
&\qquad +\int_{s}^b \left(\int_{a}^s  B_1(s,\eta) N_2(\eta,\xi)d \eta \right) x(\xi) d \xi \\
&=  \left(\mcl P_{\{0,R_{1d},R_{2d}\}} x\right)(s),
  \end{align*}}
where\vspace{\eqnspace}
\begin{align*}
R_{1d}(s,\theta)&=\int_{a}^\theta B_1(s,\xi) N_2(\xi,\theta)d \xi\\
R_{2d}(s,\theta)&=\int_{a}^s B_1(s,\xi) N_2(\xi,\theta)d \xi.\vspace{\eqnspace}
\end{align*}
The remaining two identities are minor variations on the two we most recently derived.\vspace{\eqnspace}{\scriptsize
\begin{align*}
 & \left(\mcl P_{\{0,0,B_2\}}  \mcl P_{\{0,N_1,0\}}x\right)(s)\\
&=\int_{s}^b B_2(s,\eta)\int_{a}^\eta N_1(\eta,\xi)x(\xi) d \xi d \eta\\
&=\int_{a}^b \left(\int_{a}^b \mbf I(\eta-s) \mbf I(\eta-\xi) B_2(s,\eta)  N_1(\eta,\xi)d \eta \right) x(\xi) d \xi \\
&=\int_{a}^b \mbf I(s-\xi) \left(\int_{s}^b  B_2(s,\eta) N_1(\eta,\xi)d \eta \right) x(\xi) d \xi \\
&+\int_{a}^b \mbf I(\xi-s) \left(\int_{\xi}^b  B_2(s,\eta) N_1(\eta,\xi)d \eta \right) x(\xi) d \xi \\
&=\int_{a}^s \left(\int_{s}^b  B_2(s,\eta) N_1(\eta,\xi)d \eta \right) x(\xi) d \xi \\
&+\int_{s}^b \left(\int_{\xi}^b  B_2(s,\eta) N_1(\eta,\xi)d \eta \right) x(\xi) d \xi \\
&=  \left(\mcl P_{\{0,R_{1e},R_{2e}\}} x\right)(s),
  \end{align*}}
where\vspace{\eqnspace}
\begin{align*}
R_{1e}(s,\theta)&=\int_{s}^b B_2(s,\xi) N_1(\xi,\theta)d \xi\\
R_{2e}(s,\theta)&=\int_{\theta}^b B_2(s,\xi) N_1(\xi,\theta)d \xi.\vspace{\eqnspace}
\end{align*}
Finally,\vspace{\eqnspace}{
\begin{align*}
&  \left(\mcl P_{\{0,0,B_2\}}  \mcl P_{\{0,0,N_2\}}x\right)(s)\\
&=\int_{s}^b B_2(s,\eta)\int_{\eta}^b N_2(\eta,\xi)x(\xi) d \xi d \eta\\
&=\int_{a}^b \left(\int_{a}^b \mbf I(\xi-\eta) \mbf I(\eta-s)  B_2(s,\eta)  N_2(\eta,\xi)d \eta \right) x(\xi) d \xi \\
&=\int_{a}^b \mbf I(\xi-s) \left(\int_{s}^\xi  B_2(s,\eta) N_2(\eta,\xi)d \eta \right) x(\xi) d \xi \\
&=\int_{s}^b \left(\int_{s}^\xi  B_2(s,\eta) N_2(\eta,\xi)d \eta \right) x(\xi) d \xi \\
&=  \left(\mcl P_{\{0,0,R_{2f}\}} x\right)(s),\vspace{\eqnspace}
  \end{align*}}
where\hspace{-4mm}
\[
R_{2f}(s,\theta)=\int_{s}^\theta B_2(s,\xi) N_2(\xi,\theta)d \xi.\vspace{\eqnspace}
\]
Combining all terms, we have
  \begin{align*}
  &\mcl P_{\{B_0,B_1,B_2\}}  \mcl P_{\{N_0,N_1,N_2\}}x\\
  &=\mcl P_{\{R_{0},R_{1a},R_{2a}\}}+\mcl P_{\{0,R_{1b},R_{2b}\}}+\mcl P_{\{0,R_{1c},0\}}\\
  &\quad +\mcl P_{\{0,R_{1d},R_{2d}\}}+\mcl P_{\{0,R_{1e},R_{2e}\}}+\mcl P_{\{0,0,R_{2f}\}} \\
  &=\mcl P_{\{R_{0},R_{1a}+R_{1b}+R_{1c}+R_{1d}+R_{1e},R_{2a}+R_{2b}+R_{2d}+R_{2e}+R_{2f}\}}\\
  &=\mcl P_{\{R_0,R_{1},R_{2}\}}.
  \end{align*}
The last equality holds since
\begin{align*}
&R_1(s,\theta)=R_{1a}(s,\theta)+R_{1b}(s,\theta)\\
&\qquad +R_{1c}(s,\theta)+R_{1d}(s,\theta)+R_{1e}(s,\theta)
\end{align*}
and\vspace{\eqnspace}
\begin{align*}
&R_2(s,\theta)=R_{2a}(s,\theta)+R_{2b}(s,\theta)\\
&\qquad +R_{2d}(s,\theta)+R_{2e}(s,\theta)+R_{2f}(s,\theta)
\end{align*}
$\square$
\end{pf}
This lemma proves that the class of 3-PI operators is closed under composition.\vspace{-2mm}

\begin{cor} Suppose that $\{B_i\}$ and $\{N_i\}$ are matrices of polynomials. Then if $\mcl P_{\{R_i\}}=  \mcl P_{\{B_i\}}  \mcl P_{\{N_i\}}$, $\{R_i\}$ are matrices of polynomials.\vspace{-4mm}
\end{cor}
\begin{pf}
The algebra of polynomials is closed under multiplication and integration. Therefore, the proof follows from the expressions for $\{R_i\}$ given in Lemma~\ref{thm:composition}.$\square$\vspace{-2mm}
\end{pf}

This corollary implies that the subset of 3-PI operators with polynomial parameters is likewise closed under composition and therefore forms a subalgebra.\vspace{-3mm}

\subsection{The Adjoint of 3-PI operators}\vspace{-3mm}
Next, we give a formula for the adjoint of a 3-PI operator.\vspace{-2mm}
\begin{lem}\label{lem:adjoint}
 For any bounded functions $N_0:[a,b]\rightarrow \R^{n \times n}$, $N_1,N_2:[a,b]^2 \rightarrow \R^{n \times n}$ and any $\mbf x, \mbf y \in L_2^n[a,b]$, we have\vspace{\eqnspace}
\[
\ip{\mcl P_{\{N_i\}}\mbf x}{\mbf y}_{L_2}=  \ip{\mbf x}{\mcl P_{\{\hat N_i\}}\mbf y}_{L_2}\vspace{\eqnspace}
\]
where\vspace{\eqnspace}
\begin{align}
\hat N_0(s)&=N_0(s)^T ,\quad \hat N_1(s,\eta)=N_2(\eta,s)^T, \notag \\
\hat N_2(s,\eta)&=N_1(\eta,s)^T. \label{eqn:transpose}\vspace{\eqnspace}
\end{align}
\end{lem}

\begin{pf}
Noting that\vspace{\eqnspace}
  \begin{align*}
&\ip{\mcl P_{\{N_0,N_1,N_2\}}\mbf x}{\mbf y}_{L_2}= \ip{\mcl P_{\{N_0,0,0\}}\mbf x}{\mbf y}_{L_2}\\
&\qquad +\ip{\mcl P_{\{0,N_1,0\}}\mbf x}{\mbf y}_{L_2}+\ip{\mcl P_{\{0,0,N_2\}}\mbf x}{\mbf y}_{L_2},
  \end{align*}
  we can decompose the adjoint as follows. Clearly\vspace{\eqnspace}
\[
\ip{\mcl P_{\{N_0,0,0\}}\mbf x}{\mbf y}_{L_2}=\ip{\mbf x}{\mcl P_{\{N_0,0,0\}}\mbf y}_{L_2}.\vspace{\eqnspace}
\]
For the other terms, we use a change of integration to yield\hspace{-5mm}
{\scriptsize\begin{align*}
  &\ip{\mcl P_{\{0,N_1,0\}}\mbf x}{\mbf y}=\int_{a}^b\left(\int_a^s N_1(s,\eta)\mbf x(\eta)d \eta \right)^T\mbf y(s)ds\\
  &=\int_{a}^b\int_a^b \mbf I(s-\eta) \mbf x(\eta)^T  N_1(s,\eta)^T \mbf y(s)d \eta ds\\
  &=\int_{a}^b\int_\eta^b  \mbf x(\eta)^T  N_1(s,\eta)^T \mbf y(s) ds d \eta \\
  &=\int_{a}^b\int_s^b  \mbf x(s)^T  N_1(\eta,s)^T \mbf y(\eta)  d \eta ds=\ip{\mbf x}{\mcl P_{\{0,0,\hat N_2\}}\mbf y}.\vspace{\eqnspace}
  \end{align*}}
Likewise,\vspace{\eqnspace}
 {\scriptsize \begin{align*}
  &\ip{\mcl P_{\{0,0,N_2\}}\mbf x}{\mbf y}=\int_{a}^b\left(\int_s^b N_2(s,\eta)\mbf x(\eta)d \eta \right)^T\mbf y(s)ds\\
  &=\int_{a}^b\int_a^b \mbf I(\eta-s) \mbf x(\eta)^T  N_2(s,\eta)^T \mbf y(s)d \eta ds\\
  &=\int_{a}^b\int_a^s  \mbf x(s)^T  N_2(\eta,s)^T \mbf y(\eta)  d \eta ds =\ip{\mbf x}{\mcl P_{\{0,\hat N_1,0\}}\mbf y}.\vspace{\eqnspace}
  \end{align*}}
  Combining these terms, we complete the proof. $\square$
\end{pf}

The following Corollary follows immediately from Lemma~\ref{lem:adjoint}.
\begin{cor} Suppose that $\{N_i\}$ are matrices of polynomials. Then, using the adjoint with respect to $L_2$, if $\mcl P_{\{\hat N_i\}}=  \mcl P_{\{N_i\}}^*$, $\{\hat N_i\}$ are matrices of polynomials.
\end{cor}
\vspace{-3mm}

\section{Partial Integral Equations (PIEs)}\label{sec:PIE}\vspace{-3mm}
In this section, we give the autonomous form of a Partial Integral Equation (PIE) and define notions of solution and exponential stability. Specifically, for given 3-PI operators $\mcl A,\mcl T$, we say, for an initial condition, $\mbf x_0 \in L_2^n$, that $\mbf x:\R^+ \rightarrow L_2^n$ satisfies the PIE defined by $\{\mcl A, \mcl T\}$ if $\mbf x(0)=\mbf x_0$, $\mbf x$ is Frech\'et differentiable for all $t\ge 0$ and\vspace{\eqnspace}
\begin{equation}
\mcl T \dot{\mbf x}(t)=\mcl A x(t)\label{eqn:PIE_1}\vspace{\eqnspace}\vspace{-2mm}
\end{equation}
for all $t \ge 0$.\vspace{\parspace}

Not all PIEs are well-posed in the sense of Hadamard. However, we will show in Section~\ref{sec:PDE2PIE} that if a PDE is in standardized form (satisfying Eqn.~\eqref{eqn:B_assumption}), and the PIE is generated from that standardized PDE using the formulation in Section~\ref{sec:PDE2PIE}, then the PIE is well-posed.\vspace{-2mm}

\subsection{Exponential Stability of PIEs}\vspace{-3mm}
Having defined PIE's, we now define the notion of exponential stability we will use.\vspace{-2mm}

\begin{defn}We say the PIE defined by the 3-PI operators $\{\mcl A,\mcl T\}$ is exponentially stable if there exist constants $K$ and $\gamma>0$ such that for $\mbf x(0)\in L_2^n$, any solution $\mbf x$ satisfies\vspace{\eqnspace}
\[
\norm{\mbf x(t)}_{L_2}\le K\norm{\mbf x(0)}_{L_2}e^{-\gamma t}.\vspace{\eqnspace}\vspace{-1mm}
\]
\end{defn}
In the case where the exponential stability definition holds with $\gamma=0$, we say the PIE is stable in the sense of Lyapunov or neutrally stable.\vspace{-2mm}

\section{A Unitary map from $X$ to $L_2$}\label{sec:T}\vspace{-3mm}
In this section, we show equivalence between the Hilbert space $L_2^{n_0+n_1+n_2}$ and the space\vspace{\eqnspace}
\[
X=\left\{\bmat{x_0\\x_{1}\\x_{2}}\in L_2^{n_0}\times H_1^{n_1} \times H_2^{n_2}\;:\; B {\scriptsize \bmat{ x_1(a) \\ x_1(b) \\ x_2(a) \\ x_2(b)\\ x_{2s}(a) \\ x_{2s}(b)}}=0\right\}\vspace{\eqnspace}
\]
where $B$ satisfies Equation~\eqref{eqn:B_assumption} and $X$ is equipped with the inner product\vspace{\eqnspace}
\[
\ip{\mbf x}{\mbf y}_X=\ip{x_0}{y_0}_{L_2}+\ip{\partial_s x_1}{\partial_s y_1}_{L_2}+\ip{\partial^2_s x_2}{\partial^2_s y_2}_{L_2}.\vspace{\eqnspace}
\]
Specifically, in this section, we\vspace{-2mm}
\begin{itemize}
\item Construct a unitary map $\mcl T: L_2^{n_0+n_1+n_2}\rightarrow X$ where $\mcl T$ is a 3-PI operator.
\item Show $\ip{\cdot}{\cdot}_X$ is an inner product and $X$ is Hilbert with this inner product.
\item Show that for $\mbf x \in X$, the norm $\norm{\cdot}_X$ is equivalent to the norm $\norm{\cdot}_{L_2\times H_1 \times H_2}$ where recall\vspace{\eqnspace}
\[
\norm{\mbf x}_{L_2\times H^1 \times H^2}=\norm{x_0}_{L_2}+\norm{x_1}_{H_1}+\norm{x_2}_{H_2}.\vspace{\eqnspace}
\]
\end{itemize}\vspace{-2mm}

\subsection{The Unitary Map, $\mcl T$}\vspace{-3mm}
In this subsection, we define the 3-PI operator $\mcl T=\mcl P_{\{G_i\}}$ such that if\vspace{\eqnspace}
\[
\mbf x\in X\qquad \text{and}\qquad \hat{\mbf x}\in L_2^{n_0+n_1+n_2}\vspace{\eqnspace}
\]
then\vspace{\eqnspace}
\[
\mbf x =\mcl T\bmat{I &&\\&\partial_s&\\&&\partial_s^2}\mbf x\qquad
\text{and}\qquad
\hat{\mbf x} =\bmat{I &&\\&\partial_s&\\&&\partial_s^2} \mcl T\hat{\mbf x}.\vspace{\eqnspace}
\]

\begin{figure*}[!t]
\normalsize
\begin{align*}
  G_0(s)&=\bmat{I_{n_0}&0&0\\0&0&0\\0&0&0}, & G_1(s,\theta)&=\bmat{0&0&0\\0&I_{n_1}&0\\0&0&(s-\theta)I_{n_2}}+G_2(s,\theta), & & G_2(s,\theta) =-K(s)(BT)^{-1}BQ(s,\theta),\notag\\
  G_3(s)&=\bmat{0&I_{n_1}&0\\0&0&0},\quad & G_4(s,\theta)  &=\bmat{0&0&0\\0&0&I_{n_2}}+G_5(s,\theta), & & G_5(s,\theta)=-V(BT)^{-1}BQ(s,\theta),\vspace{\eqnspace}\\[-12mm]\notag
\end{align*}
\begin{align}
&T=\bmat{I_{n_1}&0 &0\\I_{n_1}&0&0 \\0&I_{n_2} &0\\0&I_{n_2}&(b-a)I_{n_2}\\0&0&I_{n_2}\\0&0 &I_{n_2}},\quad Q(s,\theta)=\bmat{0&0&0\\0&I_{n_1} &0 \\ 0&0&0 \\0& 0&(b-\theta)I_{n_2} \\ 0& 0&0\\ 0&0& I_{n_2}},\quad K(s)=\bmat{0&0&0\\I_{n_1}&0&0\\0&I_{n_2}&(s-a)I_{n_2}}, \quad V=\bmat{0&0&0\\0&0&I_{n_2}}.\label{eqn:Gdefs}\\[-9mm]\notag
\end{align}
\hrulefill
\vspace*{4pt}
\end{figure*}

First, we first show that a PDE state $\mbf x \in H_2$ can be represented using the PIE state, $\partial_s^2\mbf x\in L_2$ and a set of `core' boundary conditions ($x(a),x_s(a)$).\vspace{-3mm}

\begin{lem}\label{lem:identity1}Suppose that $\mbf x\in H_2^n[a,b]$. Then for any $s \in [a,b]$, \vspace{\eqnspace}\vspace{-2mm}
\begin{align*}
x(s)&=x(a)+\int_a^s x_{s}(\eta)d \eta\\
x_s(s)&=x_s(a)+\int_a^s x_{ss}(\eta)d \eta\\[-1mm]
x(s)&=x(a)+x_s(a)(s-a)+\int_a^s (s-\eta)  x_{ss}(\eta)d \eta.\\[-8mm]
\end{align*}
\end{lem}
\vspace{-5mm}
\begin{pf}
  The first two identities are the fundamental theorem of calculus. The third identity is a repeated application of the fundamental theorem of calculus, combined with a change of variables. That is, for any $s \in [a,b]$,\vspace{\eqnspace}
\begin{align*}
x(s)&=x(a)+\int_a^s  x_{s}(\eta)d \eta\\
&=x(a)+\int_a^s x_s(a)ds+ \int_a^s  \int_a^\eta x_{ss}(\zeta)d \zeta d \eta.\vspace{\eqnspace}
\end{align*}
Examining the 3rd term, where recall $\mbf I(s)$ is the indicator function,\vspace{\eqnspace}
\begin{align*}
&\int_a^s  \int_a^\eta x_{ss}(\zeta)d \zeta d \eta=\int_a^b  \int_a^b \mbf I(s-\eta)\mbf I(\eta-\zeta) x_{ss}(\zeta)d \zeta d \eta\\
&=\int_a^b  \left(\int_a^b \mbf I(s-\eta)\mbf I(\eta-\zeta) d \eta\right) x_{ss}(\zeta) d \zeta\\
&=\int_a^b  \mbf I(s-\zeta)\left(\int_s^\zeta  d \eta\right) x_{ss}(\zeta) d \zeta= \int_a^s  \left(s-\zeta\right) x_{ss}(\zeta) d \zeta\vspace{\eqnspace}
\end{align*}
which is the desired result.$\square$\vspace{-2mm}
\end{pf}\vspace{-2mm}
As an obvious corollary, we have\vspace{\eqnspace}
\begin{align*}
x(b)&=x(a)+\int_a^b x_{s}(\eta)d \eta\\
x_s(b)&=x_s(a)+\int_a^b x_{ss}(\eta)d \eta\\
x(b)&=x(a)+x_s(a)(b-a)+\int_a^b (b-\eta)  x_{ss}(\eta)d \eta.\\[-7mm]
\end{align*}
The implication is that any boundary value can be expressed using two other boundary identities. In the standardized PDE representation, we have a generic set of boundary conditions defined by the matrix $B$. In the following theorem, we generalize Lemma~\ref{lem:identity1} in order to express the PDE state $\mbf x \in X$ in terms of the PIE state, $\hat{\mbf x}\in L_2^{n_0+n_1+n_2}$, and generalized BCs (which are equal to zero). This allows us to define the map $\mcl T$.\vspace{-2mm}
\begin{thm}\label{thm:identity2}
Let $\mcl T=\mcl P_{\{G_0,G_2,G_2\}}$ with $G_i$ as defined in Equations~\eqref{eqn:Gdefs}. Then for any $\mbf x \in X$, \vspace{\eqnspace}
\[
\mbf x =\mcl T \bmat{I &&\\&\partial_s&\\&&\partial_s^2} \mbf x\vspace{\eqnspace}
\]
Furthermore, for any $\hat{\mbf x},\hat{\mbf y} \in L_2^{n_0+n_1+n_2}$, $\mcl T \hat{\mbf x},\mcl T\hat{\mbf y} \in X$ and $\ip{\mcl T\hat{\mbf x}}{\mcl T \hat{\mbf y}}_X=\ip{\hat{\mbf x}}{\hat{\mbf y}}_{L_2}$.\vspace{-4mm}
\end{thm}
\begin{pf}
Suppose $\mbf x =\bmat{x_1\\x_2\\x_3}\in X$. Let us define the PIE state, $\hat{\mbf x}$, and the `full' and `core' sets of BCs as\vspace{\eqnspace}
\[
\hat{\mbf x}=\bmat{I &&\\&\partial_s&\\&&\partial_s^2} \mbf x, \quad x_{bf}={\scriptsize\bmat{ x_1(a) \\ x_1(b) \\ x_2(a) \\ x_2(b) \\ x_{2s}(a) \\ x_{2s}(b)}},\quad x_{bc}=\bmat{ x_1(a)  \\ x_2(a) \\ x_{2s}(a)}.\vspace{\eqnspace}
\]
Clearly $\hat{\mbf{x}} \in L_2^{n_0+n_1+n_2}$.
Using Lemma~\ref{lem:identity1}, we can express $x_{bf}$ using $x_{bc}$ and $\hat{\mbf{x}}$ as\vspace{\eqnspace}
\[
x_{bf}=Tx_{bc}+\mcl{P}_{\{0,Q,Q\}}\hat{\mbf{x}}\vspace{\eqnspace}
\]
Likewise, we may express $\mbf x$ in terms of $x_{bc}$ and $\hat{\mbf{x}}$ as\vspace{\eqnspace}
\[
\mbf x = K(s)x_{bc}+ \mcl{P}_{\{L_0,L_1,0\}}\hat{\mbf{x}}\vspace{\eqnspace}
\]
where\vspace{\eqnspace}
\[
L_0=\bmat{I_{n_0}&0&0\\0&0&0\\0&0&0}\qquad L_1=\bmat{0&0&0\\0&I_{n_1}&0\\0&0&(s-\theta)I_{n_2}}.\vspace{\eqnspace}
\]
We may now express the `full' boundary conditions as\vspace{\eqnspace}
\[
Bx_{bf}=BTx_{bc}+B\mcl{P}_{\{0,Q,Q\}}\hat{\mbf{x}}=0.\vspace{\eqnspace}
\]
Now from Equation~\eqref{eqn:B_assumption}, $BT$ is invertible, and hence\vspace{\eqnspace}
\begin{align*}
x_{bc}&=-(BT)^{-1}B\mcl{P}_{\{0,Q,Q\}}\hat{\mbf{x}}\\
&=-\mcl{P}_{\{(BT)^{-1}B,0,0\}}\mcl{P}_{\{0,Q,Q\}}\hat{\mbf{x}}\\
&=-\mcl{P}_{\{0,(BT)^{-1}BQ,(BT)^{-1}BQ\}}\hat{\mbf{x}}.\vspace{\eqnspace}
\end{align*}
This yields the following expression for $\mbf x$.\vspace{\eqnspace}
\begin{align*}
&\mbf x = \mcl{P}_{\{K,0,0\}}x_{bc}+ \mcl{P}_{\{L_0,L_1,0\}}\hat{\mbf{x}}\\
&\hspace{-1mm}= \hspace{-.5mm}-\mcl{P}_{\{K,0,0\}}\mcl{P}_{\{0,(BT)^{-1}BQ,(BT)^{-1}BQ\}}\hat{\mbf{x}}\hspace{-.5mm}+\hspace{-.5mm} \mcl{P}_{\{L_0,L_1,0\}}\hat{\mbf{x}}\\
&\hspace{-1mm}= \hspace{-.5mm}-\mcl{P}_{\{0,K(BT)^{-1}BQ,K(BT)^{-1}BQ\}}\hat{\mbf{x}}+ \mcl{P}_{\{L_0,L_1,0\}}\hat{\mbf{x}}\\
&\hspace{-1mm}= \hspace{-.5mm}\mcl{P}_{\{G_0,G_1,G_2\}}\hat{\mbf{x}}\vspace{\eqnspace}
\end{align*}
as desired.\vspace{\parspace}\vspace{-3mm}

Conversely, suppose that $\hat{\mbf x}=\bmat{\hat x_0 \\ \hat x_1 \\ \hat x_2} \in L_2^{n_0+n_1+n_2}$ and\vspace{\eqnspace}
\[
\mbf x=\bmat{ x_0 \\  x_1 \\  x_2} =\mcl T \bmat{\hat{ x}_0 \\ \hat{ x}_1 \\ \hat{ x}_2}.\vspace{\eqnspace}
\]
First, we note that using the definition of $\{G_i\}$, $ x_0=\hat{ x}_0 \in L_2^{n_0}$. Next, we have\vspace{\eqnspace}
\[
x_1(s)=\int_a^s \hat{x}_1(\theta)d \theta-\int_a^b \bmat{I & 0 &0}(BT)^{-1}BQ(s,\theta)\hat{\mbf x}(\theta)d \theta.\vspace{\eqnspace}
\]
Hence, $\partial_s Q(s,\theta)=0$ implies $
\partial_s  x_1 = \hat{ x}_1 \in L_2^{n_1}$.
Finally,\vspace{\eqnspace}
\begin{align*}
x_2(s)&=\int_a^s (s-\theta)\hat{ x}_2(\theta)d \theta\\
&\qquad -\int_a^b \bmat{0 & I & (s-a)}(BT)^{-1}BQ(s,\theta)\hat{\mbf x}(\theta)d \theta.\vspace{\eqnspace}
\end{align*}
Hence, since $\partial_s Q(s,\theta)=0$, $\partial_s (s-a)=1$ and $\partial_s^2 (s-a)=0$, we have $
\partial_s^2  x_2 = \hat{x}_2 \in L_2^{n_2}$.
We conclude that $\mbf x \in L_2\times H_1 \times H_2$. This now implies the identity\vspace{\eqnspace}\vspace{-1mm}
\[
Bx_{bf}=BTx_{bc}+B\mcl{P}_{\{0,Q,Q\}}\hat{\mbf x}.\vspace{\eqnspace}
\]
Using the formulae given above, we have\vspace{\eqnspace}
\[
\bmat{x_1(a)\\x_2(a)}=-\int_a^b \bmat{I & 0 &0\\0 & I & 0}(BT)^{-1}BQ(a,\theta)\hat{\mbf x}(\theta)d \theta,\vspace{\eqnspace}
\]
and\vspace{\eqnspace}
\[
x_{2s}(s)=\hspace{-1mm}\int_a^s\hspace{-1mm} \hat{x}_2(\theta)d \theta-\int_a^b \bmat{0 & 0 & I}(BT)^{-1}BQ(s,\theta)\hat{\mbf x}(\theta)d \theta,\vspace{\eqnspace}
\]
which implies\vspace{\eqnspace}\vspace{-2mm}
\[
x_{2s}(a)=-\int_a^b \bmat{0 & 0 & I}(BT)^{-1}BQ(a,\theta)\hat{\mbf x}(\theta)d \theta.\vspace{\eqnspace}
\]
Hence\vspace{\eqnspace}\vspace{-2mm}
\[
x_{bc}=\bmat{ x_1(a)  \\ x_2(a) \\  x_{2s}(a)}=-\int_a^b (BT)^{-1}BQ(a,\theta)\hat{\mbf x}(\theta)d \theta.\vspace{\eqnspace}
\]
We conclude that $BTx_{bc}=-B\mcl{P}_{\{0,Q,Q\}}\hat{\mbf x}$ and consequently\vspace{\eqnspace}\vspace{-2mm}
\[
Bx_{bf}=BTx_{bc}+B\mcl{P}_{\{0,Q,Q\}}\hat{\mbf x}=0\vspace{\eqnspace}
\]
from which we conclude that $\mbf x \in X$.\vspace{\parspace}

Finally, let $\hat{\mbf x},\hat{\mbf y} \in L_2^{n_0+n_1+n_2}$, $\mcl T \hat{\mbf x},\mcl T\hat{\mbf y} \in X$. Then\vspace{\eqnspace}
\begin{align*}
\ip{\mcl T\hat{\mbf x}}{\mcl T \hat{\mbf y}}_X
&=\ip{\bmat{I &&\\&\partial_s&\\&&\partial_s^2} \mcl T\hat{\mbf x}}{\bmat{I &&\\&\partial_s&\\&&\partial_s^2}\mcl T\hat{\mbf y}}_{L_2}\\
&=\ip{\hat{\mbf x}}{\hat{\mbf y}}_{L_2}.\\[-10mm]
\end{align*}
$\square$
\end{pf}

\begin{cor}
Let $\mcl H=\mcl P_{\{G_3,G_4,G_5\}}$ with $G_i$ as defined in Equations~\eqref{eqn:Gdefs}. Then for any $\mbf x \in X$,\vspace{\eqnspace}
\[
\bmat{0&\partial_s&0\\0&0&\partial_s}\mbf x = \mcl H\bmat{I &&\\&\partial_s&\\&&\partial_s^2}\mbf x.\vspace{\eqnspace}\vspace{-2mm}
\]
\end{cor}
\begin{pf}By Theorem~\ref{thm:identity2},\vspace{\eqnspace}
\[
\mbf x=\mcl T\bmat{I &&\\&\partial_s&\\&&\partial_s^2}\mbf x.\vspace{\eqnspace}
\]
Now, for any $\hat{ \mbf x}\in L_2^{n_0+n_1+n_2}$, it can be readily verified through differentiation that\vspace{\eqnspace}
\[
\bmat{0&\partial_s&0\\0&0&\partial_s}\mcl T\hat{\mbf x}=\mcl H \hat{\mbf x}\vspace{\eqnspace}
\]
which completes the proof. $\square$\vspace{-2mm}
\end{pf}
\begin{cor}\label{cor:unitary}
The operator $\mcl T=\mcl P_{\{G_0,G_1,G_2\}}:L_2^{n_0+n_1+n_2} \rightarrow X$ is unitary.\vspace{-2mm}
\end{cor}
\begin{pf}
Theorem~\ref{thm:identity2} shows that for any $\mbf x \in X$, there exists some $\hat{\mbf x}\in L_2$ such that $\mbf x= \mcl T \hat{\mbf x}$ (surjective). Furthermore, for any $\hat{\mbf x},\hat{\mbf y} \in L_2$, $\ip{\mcl T\hat{\mbf x}}{\mcl T \hat{\mbf y}}_X=\ip{\hat{\mbf x}}{\hat{\mbf y}}_{L_2}$, which concludes the proof.
$\square$\vspace{-2mm}
\end{pf}
Because $L_2^{n_0+n_1+n_2}$ is a Hilbert space and $\mcl T$ is unitary, Corollary~\ref{cor:unitary} implies $X$ is a Hilbert space.\vspace{-2mm}

\subsection{Equivalence of Norms}\vspace{-3mm}
In this subsection, we briefly show that the norms $\norm{\cdot}_X$ and $\norm{\cdot}_{L_2\times H_1\times H_2}$ are equivalent.

\begin{lem} For any $\mbf x \in X$, $\norm{\mbf x}_X\le \norm{\mbf x}_{L_2\times H_1\times H_2}$ and there exists a constant $c>0$ such that $\norm{\mbf x}_{L_2\times H_1\times H_2}\le c\norm{\mbf x}_X$.\vspace{-4mm}
\end{lem}
\begin{pf}
First, we note that\vspace{\eqnspace}
\begin{align*}
\norm{\mbf x}_{L_2\times H_1\times H_2}=\norm{\bmat{0\\ x_1\\ x_2}}_{L_2}+\norm{\bmat{0\\ 0\\ x_{2s}}}_{L_2}+\norm{\bmat{x_0\\ x_{1s}\\ x_{2ss}}}_{L_2}\\
=\norm{\bmat{0\\ x_1\\ x_2}}_{L_2}+\norm{\bmat{0\\ 0\\ x_{2s}}}_{L_2}+\norm{\mbf x}_{X}\\[-8mm]
\end{align*}
and hence $\norm{\mbf x}_X\le \norm{\mbf x}_{L_2\times H_1\times H_2}$. Now, since $G_i\in L_{\infty}[a,b]$, there exist $c_1,c_2>0$ such that\vspace{\eqnspace}
\begin{align*}
\norm{\bmat{0\\ x_1\\ x_2}}_{L_2}&\le \norm{\mbf x}_{L_2}=\norm{\mcl T \bmat{I &&\\&\partial_s&\\&&\partial_s^2}\mbf x}_{L_2}\\
&\le c_1\norm{\bmat{I &&\\&\partial_s&\\&&\partial_s^2}\mbf x}_{L_2}= c_1\norm{\mbf x}_{X}\\[-8mm]
\end{align*}
and\vspace{\eqnspace}
\begin{align*}
\norm{\bmat{0\\ 0\\ x_{2s}}}_{L_2}&\le \norm{\bmat{0\\ x_{1s}\\ x_{2s}}}_{L_2}=\norm{\mcl H \bmat{I &&\\&\partial_s&\\&&\partial_s^2}\mbf x}_{L_2}\\
&\le c_2\norm{\bmat{I &&\\&\partial_s&\\&&\partial_s^2}\mbf x}_{L_2}= c_2\norm{\mbf x}_{X}.\\[-8mm]
\end{align*}
Therefore, we conclude that\vspace{\eqnspace}
\[
\norm{\mbf x}_{L_2\times H_1\times H_2}\le(1+c_1+c_2)\norm{\mbf x}_{X}\vspace{\eqnspace}
\]
as desired.$\square$\vspace{-2mm}
\end{pf}
This result shows that for PDE systems in standardized form, stability in $\norm{\cdot}_X$ and $\norm{\cdot}_{L_2\times H_1\times H_2}$ are equivalent.\vspace{-2mm}

\section{Converting PDEs to PIEs}\label{sec:PDE2PIE}\vspace{-3mm}
In this section, we show that for any PDE in standardized form, there exists a PIE for which any solution to the PDE defines a solution to the PIE and any solution to the PIE defines a solution to the PDE. We further show that this result implies that exponential stability of the PIE is equivalent to exponential stability of the PDE in $X$.\vspace{-2mm}
\subsection{Equivalence of Solutions for PDEs and PIEs}\label{subsec:equivalence}\vspace{-3mm}
Now that we have the unitary 3-PI operator $\mcl T:=\mcl{P}_{\{G_0,G_1,G_2\}}$ where\vspace{\eqnspace}\vspace{-2mm}
\[
\mbf x =\mcl T \bmat{I &&\\&\partial_s&\\&&\partial_s^2} \mbf x\vspace{\eqnspace}
\]
for any $\mbf x \in X$, conversion of the PDE to a PIE (Eqn.~\eqref{eqn:PIE_1}) is direct.\vspace{-2mm}

\begin{lem} \label{lem:dynamics} Given $\hat{\mbf x}_0(t)\in L_2^{n_0+n_1+n_2}$,
the function $\hat{\mbf x}(t)\in L_2^{n_0+n_1+n_2}$ satisfies the PIE defined by $\{\mcl T,\mcl A\}$
if and only if for $\mbf x_0=\mcl T\hat{\mbf x}_0$, the function $\mbf x(t)=\mcl T \hat{\mbf x}(t)$ satisfies the PDE defined by $\{A_i,X\}$ where\vspace{\eqnspace}
\begin{align}
\mcl T&:=\mcl{P}_{\{G_0,G_1,G_2\}},\qquad \mcl A:=\mcl P_{\{H_i\}}\notag\\
H_0(s)&=A_0(s)G_0(s)+A_1(s)G_3(s)+A_{20}(s)\notag\\
H_1(s,\theta)&=A_0(s)G_1(s,\theta)+A_1(s)G_4(s,\theta),\notag\\
H_2(s,\theta)&=A_0(s)G_2(s,\theta)+A_1(s)G_5(s,\theta),\notag\\
A_{20}(s)&=\bmat{0&0&A_2(s)}\vspace{\eqnspace}\label{eqn:TAdefs}\\[-7mm]\notag
\end{align}
where the $G_i$ are as defined in Eqns.~\eqref{eqn:Gdefs}.\vspace{-2mm}
\end{lem}

\begin{pf}
Define $\mcl H:=\mcl{P}_{\{G_3,G_4,G_5\}}$. Suppose $\mbf x$ satisfies the PDE. Define\vspace{\eqnspace}
\[
\mcl D_1:=\bmat{I &&\\&\partial_s&\\&&\partial_s^2}, \quad \mcl D_2:=\bmat{0&\partial_s&0\\0&0&\partial_s}\vspace{\eqnspace}
\]
and $\hat{\mbf x}(t)=\mcl D_1 \mbf x(t)$. By Theorem~\ref{thm:identity2} and Lemma~\ref{thm:composition} and the definition of the $G_i$, we have\vspace{\eqnspace}
\begin{align*}
\dot{{\mbf x}}(t)&=\mcl{P}_{\{A_0,0,0\}}{\mbf x}(t)+\mcl{P}_{\{A_1,0,0\}}\mcl D_2{\mbf x}(t)+\mcl{P}_{\{A_{20},0,0\}}\mcl D_1 {\mbf x}(t)\\
&=\mcl{P}_{\{A_0,0,0\}}\mcl T \hat{\mbf x}(t) +\mcl{P}_{\{A_1,0,0\}}\mcl H \hat{\mbf x}(t)+\mcl{P}_{\{A_{20},0,0\}}\hat{\mbf x}(t)\\
&=\mcl{P}_{\{A_0G_0,A_0G_1,A_0G_2\}}\hat{\mbf x}(t)\\
&\quad +\mcl{P}_{\{A_1G_3,A_1G_4,A_1G_5\}}\hat{\mbf x}(t)+\mcl{P}_{\{A_{20},0,0\}}\hat{\mbf x}(t)\\
&=\mcl{P}_{\{H_0,H_1,H_2\}}\hat{\mbf x}(t)= \mcl A \hat{\mbf x}(t).\vspace{\eqnspace}
\end{align*}
Finally, $\dot{{\mbf x}}(t)=\mcl T\dot{\hat{\mbf x}}(t)$ and $\hat{\mbf x}(0)=\mcl D_1 {\mbf x}(0)=\mcl D_1 {\mbf x}_0=\mcl D_1 \mcl T{\hat{\mbf x}}_0=\hat{\mbf x}_0$.\vspace{\parspace}

Conversely, suppose $\hat{\mbf x}(t)$ solves the PIE. Define $\mbf x(t)=\mcl T \hat{\mbf x}(t)$. Then by Theorem~\ref{thm:identity2}, $\mbf x(t)\in X$ and\vspace{\eqnspace}
\begin{align*}
&\dot{{\mbf x}}(t)=\mcl T \dot{\hat{\mbf x}}(t)=\mcl A \hat{\mbf x}(t)\\
&=\mcl{P}_{\{A_0,0,0\}}\mcl T \hat{\mbf x}(t) +\mcl{P}_{\{A_1,0,0\}}\mcl H \hat{\mbf x}(t)+\mcl{P}_{\{A_{20},0,0\}}\hat{\mbf x}(t)\\
&=\mcl{P}_{\{A_0,0,0\}}{\mbf x}(t)+\mcl{P}_{\{A_1,0,0\}}\mcl D_2{\mbf x}(t)+\mcl{P}_{\{A_{20},0,0\}}\mcl D_1 {\mbf x}(t)\vspace{\eqnspace}
\end{align*}
as desired. Furthermore, ${\mbf x}(0)=\mcl T \hat{\mbf x}(0)=\mcl T \hat{\mbf x}_0={\mbf x}_0$.
$\square$\vspace{-3mm}
\end{pf}

\begin{note}
While the conversion formulae in Eqns.~\ref{eqn:Gdefs} are relatively complex, this is because they encompass a very large class of PDEs and must account for every case. Individual PIE representations of specific PDEs, by contrast are typically rather simple. In the following subsection, we demonstrate one such representation.\vspace{-2mm}
\end{note}

\subsection{PIE Representation of the E-B Beam}\label{subsec:EB_PIE}\vspace{-3mm}
To illustrate the PIE representation, we again consider the Euler-Bernoulli beam model, using the standardized PDE representation of Subsection~\ref{subsec:EB}. Applying the formulae in Eqns.~\eqref{eqn:Gdefs}, we obtain the PIE $\{\mcl T,\mcl A\}$ where\vspace{\eqnspace}
\begin{align}
\mcl T:&=\mcl P_{\{N_i\}}, & \mcl A:&=\mcl P_{\{R_i\}}&&\label{eqn:EB_PIE}\\
N_0&=0,& N_1&=\bmat{s-\theta &0\\0 &0},& N_2&=\bmat{0 & 0\\0 &\theta-s},\notag\\[-2mm]
R_0&=\bmat{0 &-c\\1 &0},& R_1&=0,& R_2&=0.\notag\\[-9mm]\notag
\end{align}

\subsection{Stability Equivalence for PDEs and PIEs}\vspace{-3mm}

\begin{lem}\label{lem:stability_equivalence}
The PDE defined by $\{A_i,X\}$ is exponentially stable in $X$ with constants $K,\gamma>0$ if and only if the PIE defined by $\{\mcl T,\mcl A\}$, where $\{\mcl T,\mcl A\}$ are as defined in Eqn.~\eqref{eqn:TAdefs}, is exponentially stable with constants $K,\gamma>0$.\vspace{-3mm}
\end{lem}
\begin{pf}
Suppose the PDE defined by $\{A_i,X\}$ is exponentially stable with constants $K,\gamma>0$. Then for any ${\mbf x}_0 \in X$, any solution ${\mbf x}$ of the PDE defined by $\{A_i,X\}$ satisfies $
\norm{{\mbf x}(t)}_{X}\le K\norm{{\mbf x}_0}_{X}e^{-\gamma t}$.
Now for $\hat{\mbf x}_0 \in L_2^{n_0+n_1+n_2}$, let $\hat{\mbf x}$ be a solution of the PIE defined by $\{\mcl T,\mcl A\}$. Define ${\mbf x}_0:=\mcl T \hat{\mbf x}_0 \in X$ and ${\mbf x}(t):=\mcl T \hat{\mbf x}(t)$. Then by Lemma~\ref{lem:dynamics}, ${\mbf x}(t)$  satisfies the PDE defined by $\{A_i,X\}$ with initial condition ${\mbf x_0}$. Therefore, by Theorem~\ref{thm:identity2},\vspace{\eqnspace}
\begin{align*}
\norm{\hat{\mbf x}(t)}_{L_2}&=\norm{\mcl T\hat{\mbf x}(t)}_{X}=\norm{{\mbf x}(t)}_{X}\\
&\le K\norm{{\mbf x}_0}_{X}e^{-\gamma t}=K\norm{\mcl T \hat{\mbf x}_0}_{X}e^{-\gamma t}\\
&=K\norm{\hat{\mbf x}_0}_{L_2}e^{-\gamma t}.\\[-8mm]
\end{align*}
Conversely, suppose the PIE defined by $\{\mcl T,\mcl A\}$ is exponentially stable with constants $K,\gamma>0$. Then for any $\hat{\mbf x}_0 \in L_2$, any solution $\hat{\mbf x}$ of the PIE defined by $\{\mcl T,\mcl A\}$ satisfies $
\norm{\hat{\mbf x}(t)}_{L_2}\le K\norm{\hat{\mbf x}_0}_{L_2}e^{-\gamma t}$. Now for $\mbf x_0 \in X$, let $\mbf x$ be a solution of the PDE defined by $\{A_i,X\}$. Define\vspace{\eqnspace}
\[
\hat{\mbf x}_0:= \bmat{I &&\\&\partial_s&\\&&\partial_s^2} \mbf x_0 \in L_2,\; \hat{\mbf x}(t):=\bmat{I &&\\&\partial_s&\\&&\partial_s^2} \mbf x(t)\in L_2.\vspace{\eqnspace}
\] Then by Lemma~\ref{lem:dynamics}, $ \mbf x(t)=\mcl T\hat{\mbf x}(t)$ and $\hat{\mbf x}(t)$  satisfies the PIE defined by $\{\mcl T,\mcl A\}$ with initial condition $\hat{\mbf x}_0$. Therefore, by Theorem~\ref{thm:identity2},\vspace{\eqnspace}
\begin{align*}
\norm{\mbf x(t)}_{X}&=\norm{\mcl T\hat{\mbf x}(t)}_{X}=\norm{\hat{\mbf x}(t)}_{L_2}\\
&\le K\norm{\hat{\mbf x}_0}_{L_2}e^{-\gamma t}=K\norm{\mcl T \hat{\mbf x}_0}_{X}e^{-\gamma t}\\
&=K\norm{\mbf x_0}_{X}e^{-\gamma t}.\\[-9mm]
\end{align*}
$\square$
\end{pf}

Having shown that PIEs are equivalent to PDEs, we now proceed to define a Linear PI Inequality (LPI), whose feasibility guarantees exponential stability a PDE in standardized form.\vspace{-2mm}

\section{Lyapunov Stability as an LPI}\label{sec:LOI}\label{sec:Stability}\vspace{-3mm}
Using the 3-PI algebra, we may now succinctly represent our Lyapunov stability conditions. The procedure is relatively straightforward.\vspace{-2mm}
\begin{thm}\label{thm:Lyapunov}
Suppose there exist $\epsilon,\delta>0$,  $N_0:[a,b]\rightarrow \R^{n \times n}$, $N_1,N_2:[a,b]^2 \rightarrow \R^{n \times n}$ such that for $\mcl P:=\mcl P_{\{N_0,N_1,N_2\}}$, $\mcl P=\mcl P^*\ge \alpha I$ and\vspace{\eqnspace}
\[
\mcl A^*\mcl P \mcl T+\mcl T^*\mcl P \mcl A \le -\delta \mcl T^*\mcl T\vspace{\eqnspace}
\]
where $\mcl T$ and $\mcl A$ are as defined in Eqn.~\eqref{eqn:TAdefs}. Then any solution, $\mbf x(t)$ of the PDE defined by $\{A_i,X\}$ satisfies\vspace{\eqnspace}
\[
\norm{\mbf x(t)}_{L_2}\le \frac{\zeta}{\alpha}\norm{\mbf x(0)}_{L_2}^2e^{-\delta/\zeta t}.\vspace{\eqnspace}
\]
where $\zeta=\norm{\mcl P}_{\mcl L(L_2)}$.
\vspace{-3mm}
\end{thm}
\begin{pf}
Suppose $\hat{\mbf x}$ solves the PIE defined by $\{\mcl T,\mcl A\}$ for some $\hat{\mbf x}_0$. Consider the candidate Lyapunov function defined as\vspace{\eqnspace}
\[
V(\hat{\mbf x})=\ip{\hat{\mbf x}(t)}{\mcl T^*\mcl P\mcl T\hat{\mbf x}(t)}_{L_2}\ge \epsilon \norm{\mcl T\hat{\mbf x}}_{L_2}^2.\vspace{\eqnspace}
\]
The derivative of $V$ along solution $\hat{\mbf x}$ is\vspace{\eqnspace}
\begin{align*}
\dot V(\hat{\mbf x}(t))&=\ip{\mcl T\dot{\hat{\mbf x}}(t)}{\mcl P \mcl T\hat{\mbf x}(t)}_{L_2}+\ip{\hat{\mbf x}(t)}{\mcl P \mcl T\dot{\hat{\mbf x}}(t)}_{L_2}\\
&=\ip{\mcl A \hat{\mbf x}(t)}{\mcl P \mcl T\hat{\mbf x}(t)}_{L_2}+\ip{\mcl T\hat{\mbf x}(t)}{\mcl P \mcl A \hat{\mbf x}(t)}_{L_2}\\
&=\ip{\hat{\mbf x}(t)}{\left(\mcl A^*\mcl P \mcl T+\mcl T^*\mcl P \mcl A\right)\hat{\mbf x}(t)}_{L_2}\\
&\le -\delta \norm{\mcl T \hat{\mbf x}(t)}_{L_2}^2.\vspace{\eqnspace}
\end{align*}
Recall $\norm{\mcl P}_{\mcl L(L_2)}=\zeta$. Then by a standard application of Gronwall-Bellman, we have\vspace{\eqnspace}
\[
\norm{\mcl T \hat{\mbf x}(t)}_{L_2}\le \frac{\zeta}{\alpha}\norm{\mcl T\hat{\mbf x}_0}_{L_2}^2e^{-\delta/\zeta t}.\vspace{\eqnspace}
\]
Now for any solution, $\mbf x$ of the PDE defined by $\{A_i,X\}$ with initial condition $\mbf x_0$, $\mbf x(t)=\mcl T\hat x(t)$ for solution of the PIE with initial condition $\hat{\mbf x}_0$ where $\mbf x_0=\mcl T\hat{\mbf x}_0$. Thus\vspace{\eqnspace}
\[
\norm{\mbf x(t)}_{L_2}\le \frac{\zeta}{\alpha}\norm{\mbf x_0}_{L_2}^2e^{-\delta/\zeta t}.\vspace{\eqnspace}\vspace{-2mm}
\]
$\square$\vspace{-2mm}
\end{pf}

\begin{note}
Theorem~\ref{thm:Lyapunov} proves exponential stability of the PDE with respect to the $L_2$ norm and not the $X$-norm. While it is possible to formulate a PIE for stability in the $X$-norm, this would differ from most existing results and the literature and hence is omitted. Note, however, that for any $\mbf x \in L_2$, $\mcl T\mbf x=0$ if and only if $\mbf x=0$ (modulo a set of zero measure).\vspace{-2mm}
\end{note}

\begin{note} Theorem~\ref{thm:Lyapunov} is equivalent to the Lyapunov inequality for PDEs with the restriction that the Lyapunov operator be a PI operator. This, in turn, may be interpreted as a dissipativity condition on the generator. Such conditions are sometimes enforced using multiplier approaches as in, e.g.~\cite{luo_book}, and have been shown to be necessary and sufficient for stability of infinite-dimensional systems, as in~\cite{datko_1970,curtain_book}. Note that the constraint that the operator $\mcl P$ be self-adjoint is not conservative as any Lyapunov function defined by a non-self-adjoint operator admits a representation using a self-adjoint operator.\vspace{-2mm}
\end{note}

Theorem~\ref{thm:Lyapunov} poses a convex optimization problem, whose feasibility implies stability of solutions of a coupled linear PDE. We refer to such optimization problems as Linear PI Inequalities (LPIs). Solving an LPI requires parameterizing the 3-PI operator $\mcl P$ using polynomials and enforcing the inequalities using LMIs. In the following section, we briefly introduce a method of enforcing positivity of a 3-PI operator using LMI constraints.\vspace{-2mm}

\section{Solving the Stability LPI via PIETOOLS}\label{sec:LPI}\vspace{-3mm}
In Section~\ref{sec:Stability}, we formulated the question of Lyapunov stability as an LPI. In this section, we will we propose a general form of LPI and show how these convex optimization problems can be solved under the assumption that all 3-PI operators are parameterized by polynomials.

For given 3-PI operators $\{\mcl E_{ij}, \mcl F_{ij},\mcl G_i\}$ and linear operator $\mcl L$, a Linear PI Inequality (LPI) is a convex optimization of the form\vspace{\eqnspace}
\begin{align}
&\min_{N_{0i}, N_{1i}, N_{2i}} \mcl L(\{N_{ij}\})\label{eqn:LPI}\\
&\sum_{j=1}^K \mcl E_{ij}^*\mcl P_{\{N_{1i},N_{2i}, N_{3i}\}} \mcl F_{ij}+\mcl G_i\ge 0 \qquad i=1,\cdots,L\notag\vspace{\eqnspace}
\end{align}
LPIs of the form of Eqn.~\eqref{eqn:LPI} can be solved directly using PIETOOLS~\cite{shivakumar_2020ACC}. Composition and adjoint are algebraic operations on the 3-PI parameters and are computed using the formulae in Section~\ref{sec:algebra}. Positivity is enforced using an LMI constraint as described in the following subsection.\vspace{-2mm}

\subsection{Enforcing Positivity of 3-PI Operators}\label{subsec:SOS}\vspace{-3mm}
In this subsection, for a given self-adjoint 3-PI operator with polynomial parameters ($\{N_i\}$), we given an LMI constraint on the coefficients of the polynomials $\{N_i\}$ which enforces a constraint of the form $\mcl P_{\{N_i\}}\ge 0$. Specifically, the following proposition (a slight modification of the result in~\cite{peet_2019TAC}) gives necessary and sufficient conditions for a 3-PI operator to have a 3-PI square root.\vspace{-2mm}

\begin{prop}\label{thm:LOI}
For any bounded functions $Z(s)$, $Z(s,\theta)$, and $g$, where $g$ is scalar and $g(s)\ge 0$ for all $s \in [a,b]$ and\vspace{\eqnspace}
{
\begin{align*}
N_0(s)&=g(s)Z(s)^T P_{11}Z(s)\\
N_1(s,\theta)&=g(s) Z(s)^T P_{12} Z(s,\theta)  + g(\theta)Z(\theta,s)^T P_{31} Z(\theta)\\
&  + \hspace{-1mm}\int_{a}^\theta g(\nu) Z(\nu,s)^T P_{33} Z(\nu,\theta) d\nu \\
& \qquad  +\hspace{-1mm}\int_\theta^s g(\nu)Z(\nu,s)^T P_{32} Z(\nu,\theta) d\nu\\
&  \qquad \qquad       +\int_s^{L}g(\nu) Z(\nu,s)^T P_{22} Z(\nu,\theta) d\nu\\
N_2(s,\theta)&=  g(s) Z(s)^T P_{13} Z(s,\theta)  + g(\theta)Z(\theta,s)^T P_{21} Z(\theta)\\
&   + \int_{a}^s g(\nu) Z(\nu,s)^T P_{33} Z(\nu,\theta) d\nu\\
&  \qquad        +\int_s^\theta g(\nu)Z(\nu,s)^T P_{23} Z(\nu,\theta) d\nu\\
&  \qquad \qquad       +\int_\theta^{L} g(\nu) Z(\nu,s)^T P_{22} Z(\nu,\theta) d\nu,\vspace{\eqnspace}
\end{align*}}
where\vspace{\eqnspace}\vspace{-2mm}
\[
P =P^T= \bmat{ P_{11} & P_{12}& P_{13}\\
      P_{21} & P_{22}& P_{23}\\
      P_{31} & P_{32}& P_{33}}\ge 0,\vspace{\eqnspace}
\]
we have $\mcl P_{\{N_i\}}^*=\mcl P_{\{N_i\}}$ and $\ip{\mbf x}{\mcl P_{\{N_i\}}\mbf x}_{L_2}\ge 0$ for all $\mbf x \in L_2[a,b]$.\vspace{-4mm}
\end{prop}
\begin{pf}
It is relatively easy to show that $\{N_i\}$ satisfy Equation~\eqref{eqn:transpose} with $\{\hat N_i\}=\{N_i\}$. Therefore, by Lemma~\ref{lem:adjoint} $\mcl P_{\{N_i\}}$ is self adjoint. Now define the operator\vspace{\eqnspace}
\[
\left(\mcl Z \mbf x\right)(s) = \bmat{\sqrt{g(s)}Z(s)\mbf x(s)\\
\int_a^s \sqrt{g(s)}Z(s,\theta)\mbf x(\theta)d\theta\\
\int_s^b \sqrt{g(s)} Z(s,\theta)\mbf x(\theta)d\theta}.\vspace{\eqnspace}
\]
Then by expanding out the composition formulae, we find $\mcl P_{\{N_i\}}=\mcl Z^* P \mcl Z$ and since $P\ge 0$, $P=(P^\half)^TP^{\half}$ for some $P^\half$. Thus \vspace{\eqnspace}
\begin{align*}
\ip{\mbf x}{\mcl P_{\{N_i\}} \mbf x}_{L_2}&=\ip{\mcl Z \mbf x}{P \mcl Z \mbf x}_{L_2}\\
&=\ip{P^{\half} \mcl Z \mbf x}{P^{\half} \mcl Z \mbf x}_{L_2}\ge 0.\\[-10mm]
\end{align*}
$\square$\vspace{-2mm}
\end{pf}
Note that Prop.~\ref{thm:LOI} does not ensure that $\mcl P_{\{N_i\}}$ is coercive. To obtain a coercive operator, one must add a coercive term of the form $\mcl P_{\{\epsilon I,0,0\}}$.\vspace{\parspace}

When we desire the $\{N_i\}$ to be polynomial, we may choose $Z$ to be the vector of monomials of bounded degree, $d$. For $g(s)=1$, the operators are positive on any domain. However, for $g(s)=(s-a)(b-s)$ the operator is only positive on the given domain $[a,b]$. For the most accurate results, we combine both choices of $g$. For notational convenience, we define the set of functions which satisfy Prop.~\ref{thm:LOI} in this way. Specifically, we denote $Z_d(x)$ as the matrix whose rows are a vector monomial basis for the vector-valued polynomials of degree $d$ or less and define the cone of positive operators with polynomial multipliers and kernels associated with degree $d$ as\vspace{\eqnspace}
\begin{align*}
&\Omega_d:=\{\mcl P_{\{N_i\}}+\mcl P_{\{M_i\}}\,:\, \label{defn:Phi} \text{ $\{N_i\}$ and $\{M_i\}$ satisfy} \notag \\
&   \text{the conditions of Prop.~\ref{thm:LOI} with $Z=Z_d$ and }\notag \\
& \text{where $g(s)=1$ and $g(s)=(s-a)(b-s)$, resp.}\}\notag\\[-7mm]
\end{align*}
The dimension of the matrices $M_i$ and $N_i$ should be clear from context. The constraint $\mcl P_{\{R_i\}} \in \Omega_d$ is then an LMI constraint on the coefficients of the polynomials $\{R_i\}$ and guarantees that $P_{\{R_i\}}\ge 0$. A Matlab toolbox (PIETOOLS) for setting up and solving LPIs based on Prop.~\ref{thm:LOI} has recently been proposed and is discussed in Subsection~\ref{sec:PIETOOLS}.\vspace{-3mm}

\subsection{The Degree-Bounded Stability Test}\label{sec:LMI}\vspace{-3mm}
By restricting the degree of the polynomial parameters, $\{N_i\}$, we obtain a PIETOOLS-based LMI which enforces the LPI conditions of Theorem~\ref{thm:Lyapunov}.\vspace{\eqnspace}
\begin{thm}\label{thm:LMI}
For any $d \in \N$, suppose there exist $\epsilon,\delta>0$, and polynomials $N_0:[a,b]\rightarrow \R^{n \times n}$, $N_1,N_2:[a,b]^2 \rightarrow \R^{n \times n}$ such that\vspace{\eqnspace}
\[
\mcl P:=\mcl P_{\{N_0-\epsilon I,N_1,N_2\}}\in \Omega_d\vspace{\eqnspace}\vspace{-2mm}
\]
and\vspace{\eqnspace}
\[
-\delta \mcl T^*\mcl T-\mcl A^*\mcl P \mcl T-\mcl T^*\mcl P \mcl A\in \Omega_d\vspace{\eqnspace}
\]
where $\mcl T$ and $\mcl A$ are as defined in Eqn.~\eqref{eqn:TAdefs}. Then any solution the PDE defined by $\{A_i,X\}$ is exponentially stable in $L_2$.\vspace{\eqnspace}
\end{thm}

Note that as mentioned in the previous subsection, the constraint $\in \Omega_d$ is an LMI constraint. \vspace{-3mm}

\subsection{PIETOOLS Implementation}\label{sec:PIETOOLS}\vspace{-3mm}
In this subsection, we give sample code using the PIETOOLS toolbox which verifies that the conditions of Theorem~\ref{thm:LMI} are satisfied.\vspace{\parspace}

A detailed manual for the PIETOOLS toolbox can be found in~\cite{shivakumar_2020ACC}. This toolbox allows for declaration and manipulation of 3-PI operators and 3-PI decision variables and enforcement of LPI constraints. PIETOOLS uses aspects of the SOSTOOLS LMI conversion process and pvar polynomial objects as defined in MULTIPOLY. PIETOOLS defines the \texttt{opvar} class of PI operators and overloads the multiplication (\texttt{*}), addition (\texttt{+}) and adjoint (\texttt{'}) operations using the formulae in Lemma~\ref{thm:composition} and Lemma~\ref{lem:adjoint}. Concatenation, and scalar multiplication are likewise defined so that 3-PI operators can be treated in a similar manner to matrices.\vspace{\parspace}

To facilitate implementation of the conditions of Theorem~\ref{thm:LMI}, we have created the script \texttt{PIETOOLS\_PDE}, which is distributed with the PIETOOLS toolbox. To use this script \textit{only} requires the user to define the standardized form of the PDE - as illustrated in Step (3) below. Specifically, the user must define \texttt{n0,n1,n2,A0,A1,A2,B}, although \texttt{A2} may be omitted if \texttt{n2=0}. The user specifies that a stability test is desired by setting \texttt{stability=1} and can specify the desired accuracy through \texttt{settings\_PIETOOLS} scripts, although by default we use \texttt{settings\_PIETOOLS\_light} script, which corresponds to $d=1$. An overview of the steps included in the \texttt{PIETOOLS\_PDE} script are provided below along with a brief description of each step.\vspace{-2mm}

	\begin{enumerate}
		\item Define independent polynomial variables. These are the spatial variables in the PDE.\vspace{-3mm}
		\begin{flalign*}
		&\texttt{pvar s,th;}&\\[-7mm]
		\end{flalign*}
		\item Initialize an optimization problem structure, $X$.\vspace{-3mm}
		\begin{flalign*}
		&\texttt{X = sosprogram([s,th]);}&\\[-7mm]
		\end{flalign*}
		\item Define the standardized PDE representation and use the provided script to construct the PIE. The script \texttt{convert\_PIETOOLS\_pde} then converts the standardized PDE to a PIE.\vspace{-3mm}
		\begin{flalign*}
		&\texttt{stability=1;}&\\
		&\texttt{n1=..;n2=..;n3=..;}&\\
		&\texttt{A0=..;A1=..;A2=..;B=..;}&\\
		&\texttt{convert\_PIETOOLS\_pde;}&\\[-7mm]
		\end{flalign*}
		\item Declare the positive 3-PI operator $\mcl P$ and add inequality constraints. $n$ is state dimension, $I$ is the interval $[a,b]$, and $d$ is the degree of the polynomial parameters in $\mcl P$. This step is encoded in the script \texttt{executive\_PIETOOLS\_stability} and is executed automatically if the user has declared the option \texttt{stability=1}.\vspace{-3mm}
		\begin{flalign*}\scriptsize
		&\hspace{-3mm}\texttt{[X,P] = poslpivar(X,n,I,d);}&\\
		&\hspace{-3mm}\texttt{D = -del*T'*T-A'*P*T-T'*P*A}&\\
		&\hspace{-3mm}\texttt{X = lpi\_ineq(X,D);}&\\[-6mm]
		\end{flalign*}
		\item Call the SDP solver. \vspace{-3mm}
		\begin{flalign*}
		&\texttt{X = sossolve(X);}&\\[-7mm]
		\end{flalign*}
		\item Get the solution. \texttt{P\_s} is the 3-PI operator, $\mcl P$.\vspace{-3mm}
		\begin{flalign*}
		&\texttt{P\_s = sosgetsol\_opvar(X,P);}&\\[-7mm]
		\end{flalign*}
	\end{enumerate}
Not that the degree, $d$, enters at Step (4) and is defined in the settings script, which defaults to \texttt{settings\_PIETOOLS\_light} ($d=1$). If higher degree is required, the setting may be changed manually or using the \texttt{settings\_PIETOOLS\_heavy} ($d=2$) script. Instructions for declaring the PDE are included in the header to \texttt{PIETOOLS\_PDE}. \vspace{-3mm}

\section{Numerical Tests of Accuracy and Scalability}\label{sec:examples1}\vspace{-3mm}
In this section, we examine the accuracy and computational complexity (scalability) of the proposed stability analysis algorithm by applying Theorem~\ref{thm:LMI} to several well-studied and relatively trivial test cases. The algorithms are implemented using the PIETOOLS toolbox described in the previous section, and use the \texttt{settings\_PIETOOLS\_light} ($d=1$) option. All computation times are listed for an Intel i7-6950x processor with 64GB RAM and only account for time taken to solve the resulting LMI using Sedumi, excluding time taken for problem setup and polynomial manipulations. In cases where the limiting value of a parameter is listed for which the system is stable, the limiting value was determined using a bisection on that parameter.\vspace{\parspace}

\noindent\textbf{Example 1:} We begin with several variations of the diffusion equation. The first is adapted from~\cite{valmorbida_2014}.\vspace{\eqnspace}
\[
\dot x(t,s)=\lambda x(t,s) + x_{ss}(t,s)\vspace{\eqnspace}
\]
where $x(0)=x(1)=0$ and which is known to be stable if and only if $\lambda <\pi^2=9.869604\cdots$. For the choice of $d=1$ in Thm.~\ref{thm:LMI}, the algorithm is able to prove stability for $\lambda\le 9.8696$ with a computation time of .54s.\vspace{\parspace}

\noindent\textbf{Example 2:} The second example from~\cite{valmorbida_2016} is the same, but changes the boundary conditions to $x(0)=0$ and $x_s(1)=0$ and is unstable for $\lambda>2.467$. For $d=1$, the algorithm is able to prove stability for $\lambda \le 2.467$ with identical computation time.\vspace{\parspace}

\noindent\textbf{Example 3:} The third example from~\cite{gahlawat_2017TAC} is not homogeneous\vspace{\eqnspace}
\begin{align*}
&\dot x(t,s)=(-.5 s^3+1.3 s^2-1.5 s+.7+\lambda) x(t,s)\\
&\qquad  + (3s^2-2s)x_{s}(t,s) + (s^3-s^2+2)x_{ss}(t,s)\\[-7mm]
\end{align*}
where $x(0)=0$ and $x_s(1)=0$ and was estimated numerically to be unstable for $\lambda > 4.65$. For $d=1$, the algorithm is able to prove stability for $\lambda \le 4.65$ with similar computation time (compare to $\lambda=4.62$ in~\cite{gahlawat_2017TAC}).\vspace{\parspace}


\noindent\textbf{Example 4:} In this example from~\cite{valmorbida_2016}, we have\vspace{\eqnspace}
\[
\dot x(t,s)=\bmat{0 & 0 & 0\\ s & 0 & 0\\ s^2 & -s^3& 0}x(t,s)+ R^{-1} x_{ss}(t,s)\vspace{\eqnspace}
\]
with $x(0)=0$ and $x_s(1)=0$. In this case, using $d=1$, we were able to prove stability for any tested value of $R$ (vs. $R\le21$ in~\cite{valmorbida_2016}) with a computation time of $4.06s$. No upper limit was found.\vspace{\parspace}

\noindent\textbf{Example 5:} For our last numerical comparison, we consider some of the recent literature on coupled linear hyperbolic systems~\cite{diagne_2012,lamare_2016,saba_2019}, often representing conservation or balance laws. Although there are several variations of the problem formulation, we consider the recent work of~\cite{saba_2019}, as it seems to contain the most accurate results. \vspace{\eqnspace}\vspace{-2mm}
\[
\dot x(t,s)=\underbrace{\bmat{0 & \sigma_1\\ \sigma_2 & 0}}_{A_0}x(t,s)+\underbrace{\bmat{-\frac{1}{r_1}&0\\0&\frac{1}{r_2}}}_{A_1}  x_{s}(t,s)\vspace{\eqnspace}
\]
with boundary conditions $x_1(0)=q x_2(0)$ and $x_2(1)=\rho x_1(1)$. In this case, we have\vspace{\eqnspace}
\[
B=\bmat{1 &-q&0&0\\0 &0&-\rho&1}.\vspace{\eqnspace}
\]
Using $d=1$, $r_1=0$, $r_2=1.1$, $\sigma_1=1$, $q=1.2$, by gridding the parameters $\sigma_2$ and $\rho$, we are able to verify stability for all stable parameter values indicated in Figure 5 in~\cite{saba_2019}. For example, at $\rho=-.4$, we were able to prove stability for $\sigma_2\le1.048$.\vspace{\parspace}

%

\noindent\textbf{Example 6 (Scalability):} Finally, we explore computational complexity using a simple $n$-dimensional diffusion equation\vspace{\eqnspace}\vspace{-2mm}
\[
\dot x(t,s)=x(t,s) +  x_{ss}(t,s)\vspace{\eqnspace}
\]
where $x(t,s) \in \R^n$. We then evaluate the computation time to perform the feasibility test for different size problems, from $n=1$ to $n=40$, choosing $d=1$ - See Fig.~\ref{fig:complexity}. Note that no factors other than $d$ influence computation time and the result is always stability.

\begin{table}\centering \begin{tabular}{c|c|c|c|c|c|c}\label{tab:taumax}
$n$ & $1$  & $5$ & $10$ & $20$ & $30$ & $40$\\
\hline
sec & .504 & 1.907 & 71.63 & 2706 & 23920 &103700  \\
\end{tabular}\caption{Number of PDEs vs. Computation Time for Stability Test}\label{fig:complexity}\vspace{-2mm}
\end{table}

\vspace{-3mm}
\section{Illustrations of Beam and Wave Equations}\label{sec:examples2}\vspace{-3mm}
In Section~\ref{sec:examples1}, we demonstrated that the proposed stability test has no obvious conservatism by finding parameter values corresponding to the stability limit for several well-studied examples. However, the representation of these PDEs in the generalized PDE for of Eqn.~\eqref{eqn:PDE_1} was straightforward. In this section, we provide some guidance on how the user might identify the $A_i$ and $B$ matrices in Eqn.~\eqref{eqn:PDE_1} for several less obvious examples - focusing on four well-known wave and beam examples. The beam examples are particularly interesting in that (to the best of our knowledge) they have not previously been analyzed using LMI-based methods. As we proceed, we call particular attention to the following two questions.\vspace{-3mm}
\begin{itemize}
  \item What are the states?
  \item What are the boundary conditions?\vspace{-3mm}
\end{itemize}
\noindent \textbf{Choice of State:} Prior to the introduction of state-space, ODEs would often be represented using scalar equations. For example, the spring-mass:\vspace{\eqnspace}
\[
m\ddot{x}(t)=-c\dot x(t)-k x(t)+F(t)\vspace{\eqnspace}
\]
is a scalar ODE. To represent this in the vector-valued state-space framework, we use $x_1$ and define an auxiliary state $x_2=\dot x$. Similarly, PDEs are often represented as scalar equations using higher-order time derivatives (e.g. The wave equation is $\ddot w=w_{xx}$).  The standardized PDE representation in Eqn.~\eqref{eqn:PDE_1}, however, uses only first-order time derivatives. Furthermore, as discussed in Subsection~\ref{subsec:EB}, the use of the standardized representation occasionally involves loss of some state information and may affect the question of stability. Specifically, the exponential stability criterion in Theorem~\ref{thm:LMI} implies all states decay exponentially. For example, If a PDE is $L_2$-stable in $u$, but not $u_s$, then if $u_s$ is included in the standardized representation, the PIE stability analysis will not be able to verify stability. 

\noindent \textbf{Boundary Conditions:} Identification of a sufficient number of boundary conditions in the universal framework is particularly important. For the $B$ matrix to have sufficient rank, the solution must be uniquely defined. One consideration to be aware of is that when we introduce additional variables to eliminate higher-order time-derivatives, these new variables must also have associated boundary conditions. This is typically solved by differentiating the original boundary conditions in time to obtain boundary conditions for the new variables.\vspace{\parspace}

In the following examples, we illustrate the process of choosing state and constructing the $A_i$ and $B$ matrices.\vspace{-3mm}

\subsection{Beam Equation Examples}\vspace{-3mm}
We first both the Euler-Bernoulli (E-B) and Timoshenko (T) models beam equations. This case is particularly interesting, as the E-B model is fundamentally diffusive and the T model has hyperbolic character. Furthermore, both these models are known to be energy-conserving\cite{luo_book}, meaning that they are stable, but not exponentially stable.\vspace{\parspace}

\noindent\textbf{Euler-Bernoulli:} In this first case, we simply recall our formulation of the cantilevered E-B beam from Subsection~\ref{subsec:EB}:\vspace{\eqnspace}
\[
\dot{\mbf x}(t)=\underbrace{\bmat{0&-c\\1&0}}_{A_2}\mbf x_{ss}(t)\vspace{\eqnspace}
\]
where $A_0=A_1=0$, $n_2=2$, and $n_0=n_1=0$. The boundary conditions take the form\vspace{\eqnspace}
\[
\underbrace{\bmat{1&0&0&0&0&0&0&0\\
0&0&0&1&0&0&0&0\\
0&0&0&0&1&0&0&0\\
0&0&0&0&0&0&0&1}}_{B}{\scriptsize\bmat{u_1(0)\\u_2(0)\\u_1(L)\\u_2(L)\\u_{1s}(0)\\u_{2s}(0)\\u_{1s}(L)\\u_{2s}(L)}}=0.\vspace{\eqnspace}
\]
Entering $\{A_i,B\}$ into the script \texttt{PIETOOLS\_PDE}, we find the E-B beam is neutrally stable (using $\delta=0$ in Thm.~\ref{thm:LMI}) for any tested value of $c>0$. However, when $\delta>0$, the code is unable to find a Lyapunov function, indicating this formulation is not exponentially stable (as expected). Note that to set $\delta=0$, we modify the script using the command \texttt{epneg=0}.\vspace{\parspace}

\noindent\textbf{Timoshenko Beam}
We now consider the Timoshenko beam model where, for simplicity, we set $\rho=E=I=\kappa=G=1$:\vspace{\eqnspace}
\begin{align*}
\ddot w &= \partial_s (w_s-\phi)&&=-\phi_s+ w_{ss}\\
 \ddot \phi &=\phi_{ss}+ (w_s-\phi)&&=  -\phi + w_s+\phi_{ss}\vspace{\eqnspace}
\end{align*}
with boundary conditions of the form\vspace{\eqnspace}
\begin{align*}
&\phi(0)=0,\qquad w(0)=0,\\
&\phi_s(L)=0,\qquad w_s(L)-\phi(L)=0.\vspace{\eqnspace}
\end{align*}
Our first step is to eliminate the second-order time-derivatives, and hence we choose $u_1=\dot w$ and $u_3=\dot \phi$. Using the boundary conditions as a guide, we choose the remaining states as $u_2=w_s-\phi$ and $u_4=\phi_s$. Note that this choice of states is a natural set of coordinates as the Timoschenko beam is known to be energy conserving with respect to these states~\cite{luo_book}. In summary, we have\vspace{\eqnspace}\vspace{-2mm}
\[
\bmat{u_1\\u_2\\u_3\\u_4}=\bmat{\dot w\\ w_s-\phi \\ \dot \phi\\ \phi_s}.\vspace{\eqnspace}
\]
This gives us 4 first order boundary conditions\vspace{\eqnspace}
\[
u_1(0)=0,\quad u_3(0)=0,\quad u_4(L)=0,\quad u_2(L)=0.\vspace{\eqnspace}
\]
Reconstructing the dynamics, we now have\vspace{\eqnspace}
\begin{align*}
  \dot{u}_{1} & = u_{2s},\quad   &&\dot{u}_{2}  =u_{1s}-u_{3}\\
  \dot{u}_{3} & = u_{4s}+ u_2,\quad   &&\dot{u}_{4}  =u_{3s}.\\[-7mm]
\end{align*}
Expressing this in our standard form we have the purely hyperbolic construction\vspace{\eqnspace}
\[
\bmat{\dot{u}_1\\ \dot u_2\\ \dot u_3\\ \dot u_4}\hspace{-2mm}=\hspace{-1mm}\underbrace{\bmat{0&0&0&0\\0&0&-1&0\\0&1&0&0\\0&0&0&0}}_{A_0}\hspace{-1mm}\bmat{u_1\\u_2\\u_3\\u_4}
+\underbrace{\bmat{0&1&0&0\\1&0&0&0\\0&0&0&1\\0&0&1&0}}_{A_1}\bmat{u_{1s}\\u_{2s}\\u_{3s}\\u_{4s}}\vspace{\eqnspace}
\]
where $A_2=[\,]$ and $n_0=n_2=0$ and $n_1=4$. The $B$ matrix is then\vspace{\eqnspace}\vspace{-2mm}
\[
\underbrace{\bmat{  1& 0 &0& 0 &0& 0& 0& 0\\
     0& 0 &1 &0 &0 &0& 0& 0\\
     0& 0& 0& 0 &0& 0 &0 &1\\
     0 &0 &0 &0 &0& 1& 0& 0}}_{B}\tiny\bmat{u_1(0)\\u_2(0)\\u_3(0)\\u_4(0)\\u_1(L)\\u_2(L)\\u_3(L)\\u_4(L)}=0\vspace{\eqnspace}
\]
where $B$ has row rank $n_1=4$ and satisfies Eqn.~\eqref{eqn:B_assumption}. The script \texttt{PIETOOLS\_PDE} indicates this system is neutrally stable (using $\delta=0$ in Thm.~\ref{thm:LMI}). However, when $\delta>0$, the code is unable to find a Lyapunov function, indicating this formulation is not exponentially stable (as expected).\vspace{\parspace}
%

\subsection{Wave Equation with Boundary Feedback Examples}\vspace{-3mm}
In this subsection, we consider wave equations attached at one end and free at the other with damping at the free end. This is a well-studied problem for which numerous stability results are available in the literature~\cite{chen_1979,datko_1986}. The simplest formulation is\vspace{\eqnspace}
\begin{align*}
  \ddot{u}(t,s) & =u_{ss}(t,s) \\
  u(t,0) & =0\qquad   u_s(t,L)=-k\dot{u}(t,L).\\[-7mm]
\end{align*}
As with the beam examples, this has a purely hyperbolic formulation. Guided by the boundary conditions, we choose\vspace{\eqnspace}
\[
  u_1(t,s)  =\dot{u}(t,s),\qquad  u_2(t,s)  =u_s(t,s).\vspace{\eqnspace}
\]This yields\vspace{\eqnspace}
\[
  \bmat{\dot u_{1}  \\ \dot u_{2}}  =\underbrace{\bmat{0&1\\1&0}}_{A_1}\bmat{u_{1s}  \\ u_{2s}}\vspace{\eqnspace}
\]
where $A_0=0$, $A_2=[\,]$ $n_1=n_2=0$ and $n_1=2$. The boundary conditions are now\vspace{\eqnspace}
\[
  \underbrace{\bmat{0&1&0&0\\0&0&k&1}}_{B}\bmat{u(0)\\ u(L)} =0.\vspace{\eqnspace}
\]
This formulation is computed to be exponentially stable in the given state $u_t,u_s$ for any tested value of $k>0$. We now consider a variation on this formulation.\vspace{\parspace}

\noindent \textbf{Diffusive Formulation} As a variation, we consider a non-diffusive formulation from~\cite{chen_1979} which was shown to be asymptotically stable in the state $u$ for $a^2+k^2>0$.\vspace{\eqnspace}
\begin{align*}
 \ddot u(t,s) &= u_{ss}(t,s)\hspace{-.5mm}-\hspace{-.5mm}2a \dot u(t,s)\hspace{-.5mm}-\hspace{-.5mm}a^2u(t,s), \, s \in [0,1] \\
  u(t,0) & =0,\qquad u_s(t,1)=-k \dot u(t,1)\vspace{\eqnspace}
\end{align*}
In this case, we are forced to choose the variables $u_1=u_t$ and $u_2=u$ yielding the diffusive formulation\vspace{\eqnspace}
\[
\bmat{\dot u_1\\ \dot u_2}=\underbrace{\bmat{-2a&-a^2\\1&0}}_{A_0}\bmat{u_1\\u_2}+\underbrace{\bmat{1\\0}}_{A_2}u_{2ss}\vspace{\eqnspace}
\]
where $A_1=0$, $n_0=0$, $n_1=1$, and $n_2=1$. Note in this case that the boundary conditions on $u_1$ force us to consider this a hyperbolic state and the boundary conditions on $u_2$ make this a diffusive state! These boundary conditions are now expressed as\vspace{\eqnspace}
\[
\bmat{0&0&1&0&0&0\\
      1&0&0&0&0&0\\
      0&k&0&0&0&1}{\scriptsize\bmat{u_1(0)\\u_1(L)\\u_2(0)\\u_2(L)\\u_{2s}(0)\\u_{2s}(L)}}=0.\vspace{\eqnspace}
\]
Computation indicates this model is neutrally stable, but not exponentially stable in the given state - a result confirmed in~\cite{chen_1979,datko_1986}.\vspace{-3mm}

\section{Conclusion}\vspace{-3mm} In this paper, we have shown how to use LMIs to accurately test stability of a large class of coupled linear PDEs. To achieve this result, we have shown how to convert well-posed coupled linear PDEs - defined on state $\mbf x_p \in X$, with associated boundary conditions and continuity constraints - to Partial-Integral Equations (PIEs) with state $\mbf x_f \in L_2$ - a formulation which is defined using the algebra of 3-PI partial-integral operators and which does not require boundary conditions or continuity constraints on $\mbf x_f$. We have shown that stability of PDEs can be reformulated as a Linear PI Inequality (LPI) expressed using 3-PI operators and operator positivity constraints. We have shown how to parameterize 3-PI operators using polynomials and how to enforce positivity of 3-PI operators using LMI constraints on the coefficients of these polynomials. We have used the Matlab toolbox PIETOOLS to solve the resulting LPIs and applied the results to a variety of numerical examples. The numerical results indicate little or no conservatism in the resulting stability conditions to several significant figures even for low polynomial degree. By conversion of LMIs developed for ODEs to LPIs, it is possible that these results can be extended to: PDEs with uncertainty; $H_{\infty}$-gain analysis of PDEs;  $H_\infty$-optimal observer synthesis for PDEs; and $H_\infty$-optimal control of PDEs. In addition, it is possible that the framework may be extended to multiple spatial dimensions using the multivariate representation proposed in~\cite{peet_2009TAC}.

\bibliographystyle{IEEEtran}
\bibliography{peet_bib,delay,PDEs}

\end{document}